\documentclass[a4paper,12pt]{amsart}

\usepackage{geometry}
\usepackage[comma,numbers]{natbib}
\usepackage[inline]{enumitem}
\usepackage{graphicx}
\usepackage[textfont=it,tableposition=above,font=small]{caption}		

\usepackage{framed}					
\usepackage[foot]{amsaddr}	
\usepackage{url}						

\usepackage{xargs}					
\usepackage{amssymb}				

\usepackage{silence}
\usepackage{todonotes}

\usepackage{pgffor}

\RequirePackage[colorlinks,citecolor=blue,urlcolor=blue]{hyperref}
\usepackage[norefs,nocites]{refcheck}
\usepackage{multirow}										
\usepackage{arydshln}										
\usepackage{psfrag}
\usepackage[symbol]{footmisc}
\usepackage{mathtools}									

\usepackage{color}

\newtheorem{theorem}{Theorem}

\theoremstyle{definition}

\newtheorem*{example}{Example}
\newtheorem*{remark}{Remark}

\newcommand{\Prob}[1][\R^d]{\mathcal P\left(#1\right)}
\newcommand{\PP}{\mathsf P}														
\newcommand{\as}{\mathrm {a.s.}}											
\newcommand{\R}{\mathbb R}														
\newcommand{\dd}{\,\mathrm{d}\,}											
\newcommand{\tr}{^\mathsf{T}}													
\newcommand{\D}{HD}																		
\newcommand{\PD}{PD}																	
\newcommand{\Sph}[1][d]{\mathbb{S}^{#1-1}}						
\newcommand{\SphO}{\mathbb{S}^{1}}										

\title[Convergence rates for the approximated depth]{Uniform convergence rates for the approximated halfspace and projection depth}

\author{Stanislav Nagy$^1$}

\address{\hspace{-1em}$^1$Charles University,
	Faculty of Mathematics and Physics,
	Department of Probability and Math. Statistics,
	Prague, Czech Republic
}

\email{nagy@karlin.mff.cuni.cz}

\author{Rainer Dyckerhoff$^{2}$}

\address{$^2$University of Cologne, 
				 Institute of Econometrics and Statistics, 
         Cologne, Germany}

\author{Pavlo Mozharovskyi$^{3}$}

\address{$^3$LTCI, T\'el\'ecom Paris, 
				 Institut polytechnique de Paris, France}

\date{\today}

\begin{document}

\begin{abstract}
The computational complexity of some depths that satisfy the projection property, such as the halfspace depth or the projection depth, is known to be high, especially for data of higher dimensionality. In such scenarios, the exact depth is frequently approximated using a randomized approach: The data are projected into a finite number of directions uniformly distributed on the unit sphere, and the minimal depth of these univariate projections is used to approximate the true depth. We provide a theoretical background for this approximation procedure. Several uniform consistency results are established, and the corresponding uniform convergence rates are provided. For elliptically symmetric distributions and the halfspace depth it is shown that the obtained uniform convergence rates are sharp. In particular, guidelines for the choice of the number of random projections in order to achieve a given precision of the depths are stated. 
\end{abstract}

\maketitle

\section{Introduction}

Data depth is a general concept that intends to provide a base for nonparametric statistical analysis of multivariate data. The idea is to quantify the centrality of each point $x$ in $\R^d$, with respect to a given dataset. Points that tend to lie in the central bulk of the data are assigned high depth values; points that fail to follow the prevailing pattern of the observations are flagged by their low depth. One obtains a data-dependent ranking of points in $\R^d$, that enables constructions of analogues of quantiles, ranks, and orderings, applicable to multivariate datasets. Many definitions of depths appeared in the literature; we refer to \cite{Zuo_Serfling2000, Liu_etal1999, Mosler2002} and references therein. Here we focus on two important depths that share common traits --- the halfspace depth of \citet{Tukey1975}, and the (generalized) projection depth whose history traces back to \citet{Stahel1981} and \citet{Donoho1982}.

While theoretically the depths are appealing, it is in practice often difficult to evaluate them exactly. For instance, the computation of the halfspace depth of a single point in arbitrary dimension is known to be NP-hard \cite{Johnson_Prepata1978}. Therefore, a great deal of research has focused on procedures that approximate the true depth \cite{Dyckerhoff2004, Cuesta_Nieto2008, Chen_etal2013, Shao_Zuo2012, Bogicevic_Merkle2018, Shao_Zuo2019}. A particularly simple upper bound on the halfspace and projection depth can be devised if one uses their so-called \emph{projection property} \citep{Dyckerhoff2004}, which means that the overall (multivariate) depth of a point $x$ is expressed as the infimum of (univariate) depths of projections of $x$ with respect to the projected dataset. This suggests the following approximation procedure: \begin{enumerate*}[label=(\roman*)]\item draw a random sample of $n$ directions $U_i$, $i=1,\dots, n$, uniformly from the unit sphere $\Sph$ of $\R^d$; \item evaluate the (univariate) depths of $\left\langle U_i, x \right\rangle$ with respect to the dataset projected onto $U_i$ for each $i$; and \item approximate the depth of $x$ by the minimum of these numbers.\end{enumerate*} This approximation was first proposed in \citep{Dyckerhoff2004}; for the halfspace depth, it leads to what is sometimes called the random Tukey (or halfspace) depth \citep{Cuesta_Nieto2008}.

Here we study the statistical properties of this approximation procedure. We are interested in conditions under which uniform convergence of the approximated depth to its theoretical counterpart can be guaranteed, and its uniform convergence rates. It turns out that for finite datasets, uniform convergence is never achieved. Thus, we focus on general probability measures on $\R^d$, and show that under appropriate regularity conditions, uniform approximations of the depths are valid, and their sharp convergence rates are possible to be devised. These results lend valuable insights into the behaviour of the approximated depths. They allow us to state general guidelines for the number of directions $n$ that need to be taken in order to achieve a desired precision. Nevertheless, even in very regular models, $n$ grows fast with the dimension $d$. Above $d=5$, for reasonable $n$ the considered randomization scheme does not attain precision sufficient for practical applications, and more elaborate approximation methods should be preferred.

In Sections~\ref{section:halfspace depth}, \ref{section:main results} and \ref{section:examples} we deal with the halfspace depth. Section~\ref{section:halfspace depth} introduces notations and fixes the basic ideas of the paper. Our two main theorems, of a rather technical nature, are given in Section~\ref{section:main results}. Their applications to different classes of probability distributions are discussed in Section~\ref{section:examples}. Explicit, and exact, rates of convergence are derived for \begin{enumerate*}[label=(\roman*)] \item general unimodal elliptically symmetric distributions; \item multivariate Gaussian distributions; \item multivariate $t$-distributions; \item uniform distributions on balls; and \item the general collection of multivariate $p$-symmetric measures \citep{Fang_etal1990}. \end{enumerate*} We give explicit guidelines for the choice of the parameter $n$ that allow to achieve a pre-specified quality of the approximation. In Section~\ref{section:discontinuous case} we discuss situations when uniform approximation cannot be achieved; in particular, we deal with empirical measures. Extensions to the (generalized, or asymmetric) projection depth are treated in Section~\ref{section:projection depth}. Some concluding remarks can be found in Section~\ref{section:other sampling distributions}. All the proofs are deferred to the appendix --- Appendix~\ref{appendix:halfspace depth} for the halfspace depth, and Appendix~\ref{appendix:projection depth} for the generalized projection depths.

\section{Halfspace depth and its approximations}	\label{section:halfspace depth}

Let $\left(\Omega, \mathcal A, \PP\right)$ be the probability space on which all random elements are defined. Denote by $\Prob$ the space of all Borel probability measures on $\R^d$ equipped with the Euclidean norm $\left\Vert \cdot \right\Vert$ and the scalar product $\left\langle \cdot, \cdot \right\rangle$, $d = 1, 2, \dots$. We write $\Sph = \left\{ x \in \R^d \colon \left\Vert x \right\Vert = 1 \right\}$ for the unit sphere in $\R^d$. Notation $X \sim P$ stands for a random variable $X$ whose distribution is $P \in \Prob$.

For $P \in \Prob$ and $x \in \R^d$, the \emph{halfspace depth} (or \emph{Tukey depth}) of $x$ with respect to $X \sim P$ is defined as
	\begin{equation}	\label{halfspace depth}
	\D(x;P) = \inf_{u \in \Sph} \PP\left( \left\langle u, X \right\rangle \leq \left\langle u, x \right\rangle \right).	
	\end{equation}
Where no mistake can be made, $\D(x;P)$ is shortened to $\D(x)$. The halfspace depth was proposed in \cite{Tukey1975}, see also \cite{Donoho1982, Donoho_Gasko1992} and \cite{Rousseeuw_Ruts1999}.

For $n = 1,2,\dots$ consider $U_1, \dots, U_n$ a random sample from the uniform distribution on $\Sph$. We are interested in the quality of the approximation of the depth \eqref{halfspace depth} by its randomized counterpart \cite[Section~6]{Dyckerhoff2004}
	\begin{equation}	\label{approximated depth}
	\D_n(x;P) = \min_{i=1,\dots,n} \PP\left( \left\langle U_i, X \right\rangle \leq \left\langle U_i, x \right\rangle \right).
	\end{equation}
Again, $\D_n(x;P)$ will be shortened to $\D_n(x)$. It is evident that
	\[	\D_n(x) \geq \D_{n+1}(x) \geq \D(x) \quad \mbox{for all $x\in\R^d$ and $n=1,2\dots$},	\]
and that for $d=1$ the random depth $\D_n$ with high probability reduces to the true depth $\D$
	\begin{equation*}	
	\PP\left( \sup_{x\in\R} \left\vert \D_n(x) - \D(x) \right\vert = 0 \right) = 1 - 2 \left(\frac{1}{2}\right)^n \quad \mbox{for all $n=1,2,\dots$}	
	\end{equation*}
In the interesting case $d>1$, for any $x\in\R^d$ the random $\D_n(x)$ approximates $\D(x)$. In \cite[Proposition~11]{Dyckerhoff2004} it was shown that for any $x \in \R^d$ and $P \in \Prob$ the convergence $\D_n(x) \xrightarrow[n\to\infty]{\as} \D(x)$ holds true. Here we are interested in establishing uniform extensions of that convergence result, and deriving the corresponding rates of convergence; for an array of applications of these results to the computation of the depth see \cite[Section~2.3]{Mozharovskyi_etal2015}.

Define the \emph{halfspace function} of $X \sim P$, given for $x \in \R^d$ by
	\begin{equation}	\label{phi function}
	\varphi_x \colon \Sph \to [0,1] \colon u \mapsto \PP\left( \left\langle u, X \right\rangle \leq \left\langle u, x \right\rangle \right).
	\end{equation}
Both $\D(x)$ and its approximation $\D_n(x)$ can be expressed in terms of $\varphi_x$ as
	\begin{equation*}	
	\begin{aligned}
	\D(x) = \inf_{u \in \Sph} \varphi_x(u), \quad \mbox{and} \quad \D_n(x) = \min_{i=1,\dots,n} \varphi_x(U_i).
	\end{aligned}
	\end{equation*}
If the minimum value of $\varphi_x$ is attained in a single direction in $\Sph$, we denote this direction by $\widetilde{u}(x) \in \Sph$, that is $\D(x) = \varphi_x\left(\widetilde{u}(x)\right)$.


When deriving the rates of convergence of the depth approximations that are uniform on a set $C\subset\R^d$, it is necessary to assume a certain form of the equicontinuity of the halfspace functions
	\begin{equation}	\label{phi functions}
	\left\{ \varphi_x \colon x \in C \right\}.
	\end{equation}
In the first step, in Section~\ref{section:Lipschitz continuity} we derive our main result under the assumption of uniformly Lipschitz continuous functions $\varphi_x$. In the second step in Section~\ref{general continuity} we extend the previous theorem, and deal with a general equicontinuous class of functions \eqref{phi functions}.

\section{Approximation of the halfspace depth: Main results}	\label{section:main results}

We are interested in the rate of convergence of the random sequence
	\[	\Delta_n(C) = \sup_{x \in C} \left\vert \D_n(x) - \D(x) \right\vert	\]
as $n\to\infty$. We start by two technical results that will be of great importance in the sequel. For $\Gamma\left(\cdot\right)$ the gamma function, denote by
	\begin{equation}\label{a_d function}
        a_d(\varphi) = \frac{\Gamma\left(d/2\right)}{\Gamma\left((d-1)/2\right) \sqrt{\pi}} \int_0^\varphi \left(\sin(\theta)\right)^{d-2} \dd \theta, \quad \mbox{for }\varphi \in [0,\pi]	
        \end{equation}
the $(d-1)$-dimensional Hausdorff measure of a cap of a polar angle\footnote{The angle between the rays from the center of the hypersphere to the apex of the cap, and the edge of the $(d-1)$-dimensional disk that forms the base of the cap.} $\varphi$ of a sphere in $\R^d$ with unit surface (hyper-)area. In other words, $a_d\left(\varphi\right)$ is the ratio of the surface area of the spherical cap of polar angle $\varphi$ in $\Sph$, and the surface area of the whole $(d-1)$-sphere $\Sph$. The function $a_d$ takes particularly simple forms for $d = 2$ and $d = 3$
	\[	a_2\left(\varphi\right) = \frac{\varphi}{\pi} \quad \mbox{and} \quad a_3\left(\varphi\right) = \left(\sin\left(\frac{\varphi}{2}\right)\right)^2.	\]
For higher dimensions $d$, it is a trigonometric polynomial. Note that all functions $a_d$ are strictly increasing and continuous on their domain. Denote the inverse of $a_d$ by $a_d^{-1} \colon [0,1] \to [0,\pi]$. 

\subsection{Lipschitz continuous halfspace functions} 	\label{section:Lipschitz continuity}

For starters, assume that the set of functions \eqref{phi functions} is uniformly Lipschitz continuous with constant $L \geq 0$, i.e. that
	\begin{equation}	\label{Lipschitz continuity}
	\sup_{x\in C} \left\vert \varphi_x(u) - \varphi_x(v) \right\vert \leq L \left\Vert u - v \right\Vert_g \quad \mbox{for all }u, v \in \Sph,	
	\end{equation}
where $\left\Vert u - v \right\Vert_g = \arccos\left(\left\langle u, v \right\rangle\right)$ is the great-circle distance, i.e. the geodesic distance between $u$ and $v$. Note that the great-circle distance and the Euclidean distance are related via
	\[
	\left\Vert u - v \right\Vert
	=\sqrt{2\left(1-\cos\left(\left\Vert u - v \right\Vert_g\right)\right)}
	=2\sin\left(\left\Vert u - v \right\Vert_g/2\right)\,.
	\]
This relation follows easily from the cosine formula and the half-angle formula. A condition similar to \eqref{Lipschitz continuity} was used in \cite{Burr_Fabrizio2017} when deriving the uniform rates of convergence for the ordinary empirical halfspace depth process. Examples of distributions that satisfy property \eqref{Lipschitz continuity} will be listed in Section~\ref{section:examples}. The proof of the following result can be found in Appendix~\ref{appendix:Lipschitz}.

\begin{theorem}	\label{theorem:Lipschitz}
Let $P\in\Prob$ be such that \eqref{Lipschitz continuity} holds true for a set $C \subset \R^d$. Then
	\[	{\lim\sup}_{n\to\infty} \frac{n \, a_d \left(\Delta_n(C)/L\right) - \log n}{\log \log n} \leq d \quad \as	\]
\end{theorem}

\subsection{General uniformly continuous projections}	\label{general continuity}

Now, assume that the class of functions \eqref{phi functions} is uniformly equicontinuous, but not necessarily uniformly Lipschitz. Consider the minimal modulus of continuity of this class of functions 
	\begin{equation} 	\label{modulus}
	\delta \colon [0,\pi] \to [0,1] \colon t \mapsto \sup_{x\in C} \sup_{u, v \in \Sph, \, \left\Vert u - v \right\Vert_g \leq t} \left\vert \varphi_x(u) - \varphi_x(v) \right\vert.	
	\end{equation}
Using the minimal modulus of continuity we can bound the distance between
$\varphi_x(u)$ and $\varphi_x(v)$ uniformly for all $x\in C$,
	\begin{equation}	\label{modulus bound}
	\sup_{x\in C} \left\vert \varphi_x(u) - \varphi_x(v) \right\vert \leq \delta\left(\left\Vert u - v \right\Vert_g\right) \quad \mbox{for all }u, v \in \Sph.	
	\end{equation}
Under the condition of uniform equicontinuity of \eqref{phi functions}, the function $\delta(t)$ is continuous, non-negative and non-decreasing, with $\delta(0) = 0$, see \cite[Chapter~2, \S~6]{deVore_Lorentz1993}. With no loss of generality, in the sequel we assume that $\delta$ is increasing, and denote its inverse function by $\delta^{-1}$. If $\delta$ is constant on some interval, proper use of a generalized inverse function of $\delta$ in place of $\delta^{-1}$ yields the same conclusions \citep{Embrechts_Hofert2013}. The proof of the next result is in Appendix~\ref{appendix:modulus}.

\begin{theorem}	\label{theorem:modulus}
Let $P\in\Prob$ and $C \subset \R^d$ be such that the function $\delta$ from \eqref{modulus} is continuous at $0$ with an inverse function $\delta^{-1}$. Then
	\[	
	{\lim\sup}_{n\to\infty} \frac{n \, a_d \left(\delta^{-1}\left(\Delta_n(C)\right)\right) - \log n}{\log \log n} \leq d \quad\as	
	\] 
\end{theorem}

If \eqref{phi functions} is uniformly H\"older continuous, i.e. for some $K>0$ and $\alpha \in (0,1]$ we may choose $\delta(t) \leq K t^\alpha$, the last formula gives
	\[	{\lim\sup}_{n\to\infty}  \frac{n \, a_d \left( \left(\Delta_n(C) / K\right)^{1/\alpha} \right)  - \log n}{\log\log n} \leq d \quad\as	\]

\section{Applications to probability distributions}	\label{section:examples}

\subsection{Approximations and affine invariance}	\label{section:affine invariance}

A depth $D$ is said to be affine invariant if for any $X \sim P \equiv P_X \in \Prob$, $A \in \R^{d \times d}$ non-singular, and $b \in \R^d$ 
	\[	D(x;P_X) = D(A x + b; P_{A X + b}) \qquad \mbox{for all }x\in\R^d.	\]
Affine invariance is a desired property of a depth function. Both the halfspace depth and the (generalized) projection depth satisfy it. 

The following theorem asserts that for the halfspace depth, computation of approximations with respect to affine images of a measure affects the resulting rates of convergence only by a multiplicative constant. This allows us to focus in some examples in the sequel only on the isotropic situation, which typically means that the measure is considered to be centred to have the expected value in the origin, and the scatter matrix a multiple of the identity matrix. The proof of the next result can be found in Appendix~\ref{appendix:affine}.

\begin{theorem}	\label{theorem:affine}
Let Theorem~\ref{theorem:modulus} hold true. Then for $A \in \R^{d \times d}$ non-singular, $b\in\R^d$, and $X \sim P \in \Prob$ there exists a constant $K > 0$ that depends only on $A$ such that for
	\[	\widetilde{\Delta}_n(C) = \sup_{x \in A C + b} \left\vert \D_n(x;P_{A X + b}) - \D(x;P_{A X + b}) \right\vert	\]
it holds true that
	\[
	{\lim\sup}_{n\to\infty} \frac{n \, a_d \left(\delta^{-1}\left(\widetilde{\Delta}_n(C)\right)/K\right) - \log n}{\log \log n} \leq d \quad\as	
	\]
\end{theorem}

\subsection{General distributions with densities}	

Let us now explore when the equicontinuity conditions \eqref{Lipschitz continuity} or \eqref{modulus} hold true. The first interesting case is that of a bounded set $C \subset \R^d$.

\begin{theorem}	\label{theorem:density example}
Let $P\in\Prob$ be a distribution such that
	\begin{equation}	\label{Delta}
	P(u\tr X = \alpha) = 0 \mbox{ for all $u\in\Sph$ and $\alpha \in \R$}
	\end{equation}	
that is concentrated in a bounded subset $C$ of $\R^d$, i.e. $P(C) = 1$. Then $\Delta_n\left(\R^d\right) \xrightarrow[n\to\infty]{\as}0$.
\end{theorem}

Condition \eqref{Delta} is satisfied if $P$ admits a density in $\R^d$. If the density is bounded, an explicit rate of convergence of $\Delta_n(C)$ can be provided.

\begin{theorem}	\label{theorem:bounded example}
Let $P\in\Prob$ be a distribution with a density bounded from above by a constant $M>0$, that vanishes outside a bounded set $C \subset \R^d$. Let further be $\mathrm{diam}(C)=\sup\{\left\Vert x-y \right\Vert \colon x,y\in C\}$ the diameter of $C$. Then the rate from Theorem~\ref{theorem:Lipschitz} holds true with 
	\[	L = \frac{M \pi^{d/2-1} \left(\mathrm{diam}(C)\right)^d}{\Gamma\left(d/2+1\right)}.	\]
\end{theorem}

For the proofs of these two theorems see Appendix~\ref{appendix:density example} and~\ref{appendix:bounded example}. Note that the rate from Theorem~\ref{theorem:bounded example} is quite general, yet usually weak.

\subsection{Elliptically symmetric distributions}	\label{section:elliptically symmetric distributions}

A much finer result can be shown for elliptically symmetric distributions with unimodal densities on $\R^d$. By a unimodal elliptically symmetric density we understand a density
	\[	f(x) = \left\vert \Sigma \right\vert^{-1/2} g\left( \left(x - \mu\right)\tr \Sigma^{-1} \left( x - \mu \right) \right) \quad\mbox{for }x \in \R^d,	\]
with $\Sigma \in \R^{d \times d}$ non-singular, $\mu \in \R^d$, and $g \colon [0,\infty) \to [0,\infty)$ a non-increasing scalar function. This important collection of distributions covers all Gaussian distributions, or multivariate extensions of the $t$-distributions. In this situation we are allowed to consider $C = \R^d$, and study the uniform convergence of the approximated depth over the whole sample space. In accordance with the discussion in Section~\ref{section:affine invariance} we primarily focus on distributions with $\Sigma$ a multiple of the identity matrix. Such distributions are frequently called spherically symmetric distributions.

\begin{theorem}	\label{theorem:elliptical}
Let $P\in\Prob$ be a distribution with an elliptically symmetric density and $\Sigma$ a multiple of the identity matrix. Then $\Delta_n(\R^d) \xrightarrow[n\to\infty]{\as}0$. If, moreover, the density of $P$ is unimodal, then the rate in Theorem~\ref{theorem:modulus} holds true with both function $\delta$ taken as $\delta_1(t) = (1-\cos(t))/2$, and $\delta_2(t) = t^2/4$, with $C = \R^d$.
\end{theorem}

Note that in the proof of the previous theorem, given in Appendix~\ref{appendix:elliptical}, for $F_p$ the distribution function of the random variable $\left\langle u, X \right\rangle$, $u \in \Sph$, the better bound 
	\[	\sup_{t \geq 0} \left( F_p(t) - F_p(t\cos(\varepsilon)) \right) \leq (1-\cos(\varepsilon))/2	\]
is still not the tightest possible one. Thus, Theorem~\ref{theorem:elliptical} is suboptimal for particular choices of $P$. For elliptically symmetric distributions, define the \emph{tight modulus of continuity} of the halfspace functions (the modulus of continuity of $F_p$ evaluated near the argument of the minimum of $\varphi_x$) by
	\begin{equation}	\label{delta}
	\delta(\varepsilon) = \sup_{t \geq 0} \left( F_p(t) - F_p(t\cos(\varepsilon)) \right) \quad \mbox{for }\varepsilon \in [0,\pi].
	\end{equation}
Using the same proof as before with this choice of the modulus $\delta$, it can be asserted that for all $n$ large enough the rate from Theorem~\ref{theorem:modulus} holds true for $P$ also with $\delta$ given in \eqref{delta} and $C = \R^d$.

Interestingly, the latter bound using the modulus \eqref{delta} can be shown to be optimal for spherically symmetric distributions. For this, recall that for two sequences of real numbers $a = \left\{a(n)\right\}_{n=1}^\infty$ and $b = \left\{b(n)\right\}_{n=1}^\infty$ we say that $a$ is asymptotically bounded both from above and below by $b$, written $a(n) = \Theta(b(n))$, if 
	\[	0 < {\lim\inf}_{n\to\infty} \frac{a(n)}{b(n)} \leq {\lim\sup}_{n\to\infty} \frac{a(n)}{b(n)} < \infty.	\]
The proof of the next theorem is given in Appendix~\ref{appendix:optimal}.

\begin{theorem}	\label{theorem:optimal}
Let $P\in\Prob$ be a distribution with a unimodal elliptically symmetric density with $\Sigma$ a multiple of the identity matrix, and let $\delta(\varepsilon)$ be defined as in \eqref{delta}. 
Then
	\[	{\lim\inf}_{n\to\infty} \frac{n \, a_d \left(\delta^{-1}\left(\Delta_n(\R^d)\right)\right) - \log n}{\log \log n} \geq d - 2 \quad \as	\] 
In particular, for $d>2$
	\[	\PP\left( n \, a_d \left(\delta^{-1}\left(\Delta_n(\R^d)\right)\right) - \log n = \Theta \left( \log\log n \right) \right) = 1.	\]
\end{theorem}

Specific rates of convergence can be explored with particular models at hand. Here we study three important special cases --- multivariate Gaussian distributions, multivariate elliptically symmetric $t$-distributions, and uniform distributions on ellipsoids. In accordance with the discussion from Section~\ref{section:affine invariance}, we treat only standardized distributions.

\begin{figure}[htpb]
	\includegraphics[width=.45\textwidth]{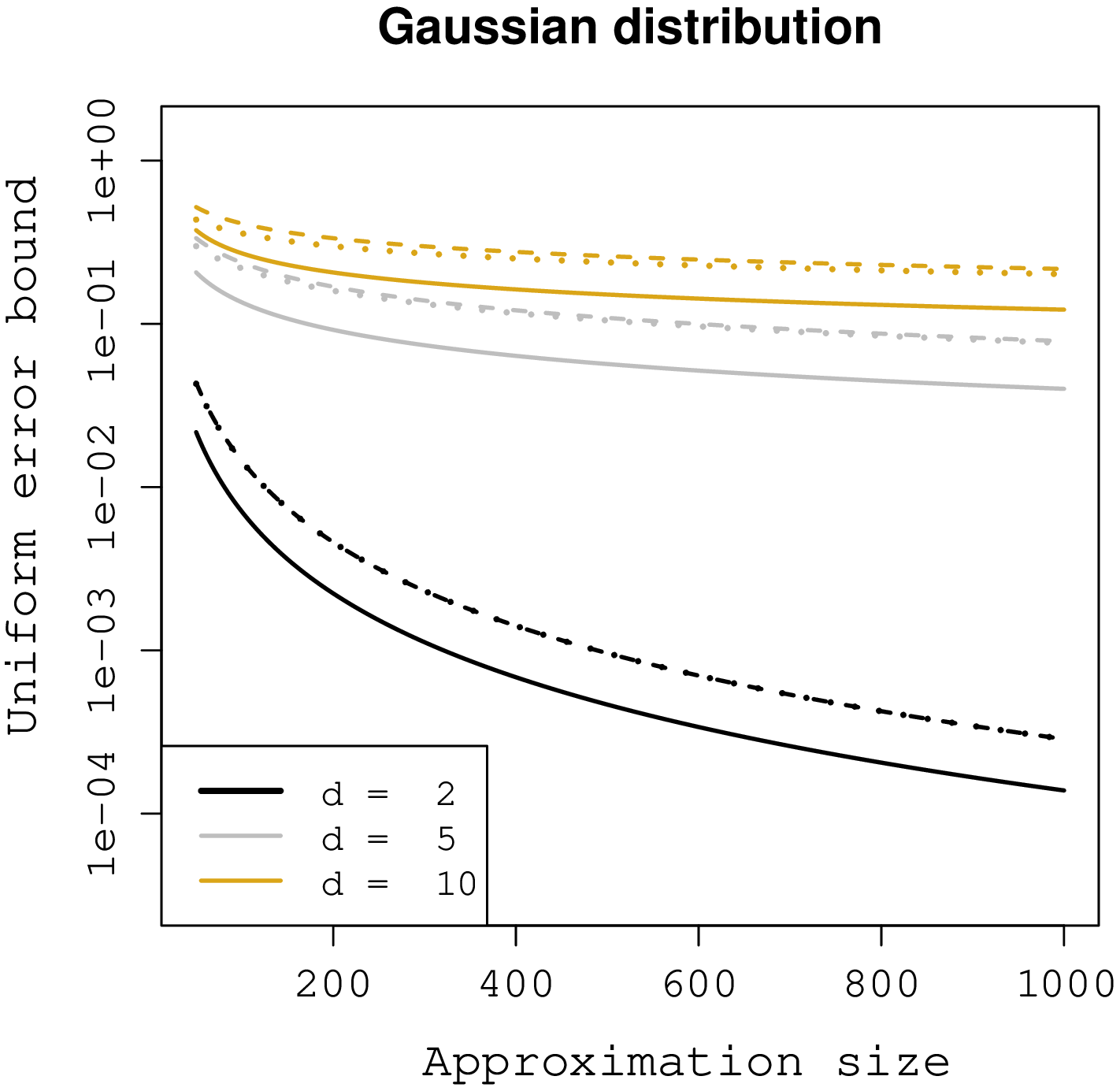} \qquad \includegraphics[width=.45\textwidth]{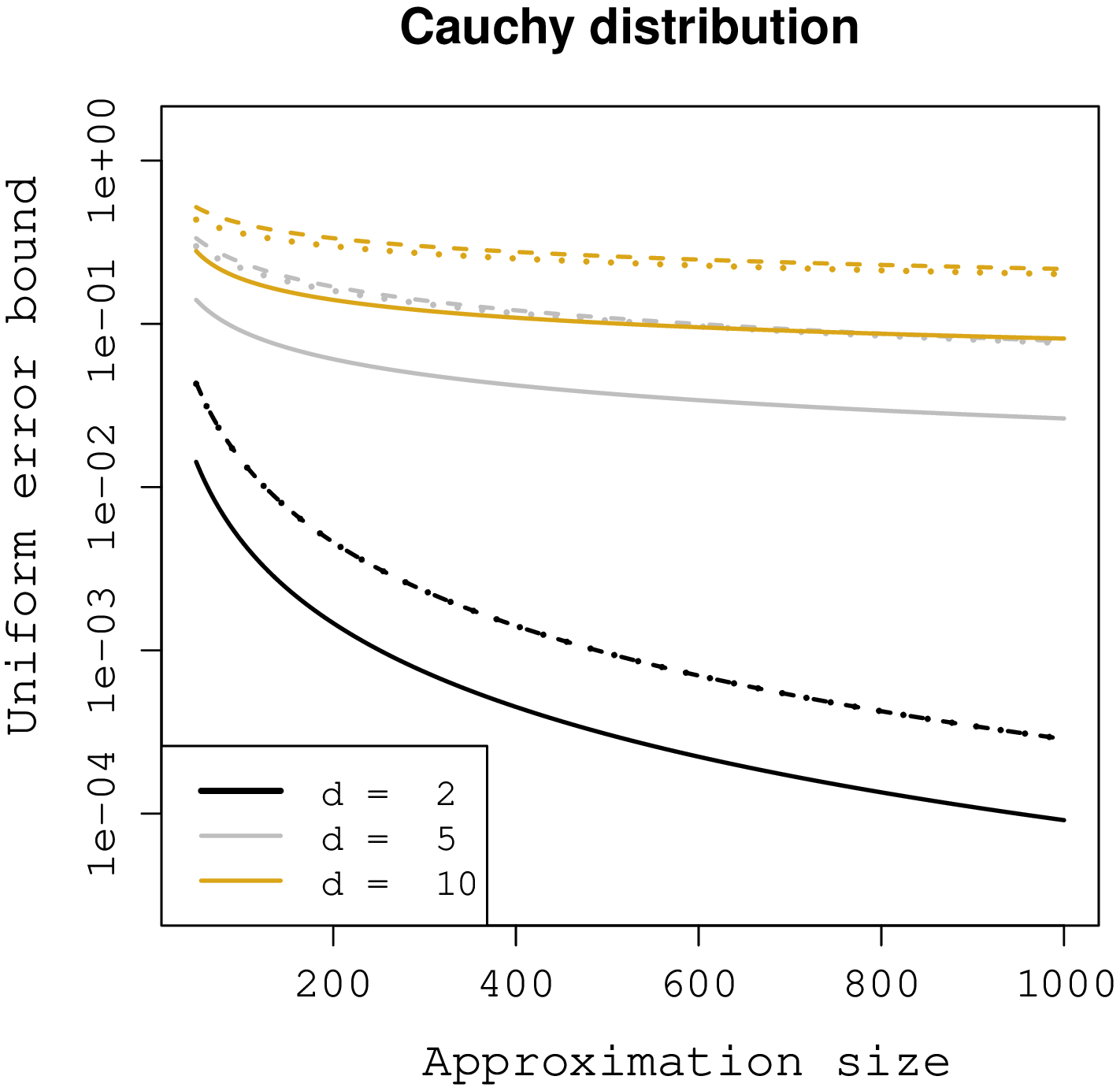} 
	\caption{Uniform error bounds for $\Delta_n(\R^d)$ as a function of $n$: the (standard) multivariate Gaussian distribution (left panel) and the (standard) multivariate Cauchy distribution (right panel). The solid curves represent the tight error bounds from Theorem~\ref{theorem:optimal}; the dashed and the dotted lines represent the two respective bounds from Theorem~\ref{theorem:elliptical}, valid for any unimodal elliptically symmetric distribution. The dashed and the dotted lines are nearly identical in the plot. Note the logarithmic scale of the vertical axes.}
	\label{figure:Gaussian_Cauchy}
\end{figure}

\begin{figure}[htpb]
	\includegraphics[width=.65\textwidth]{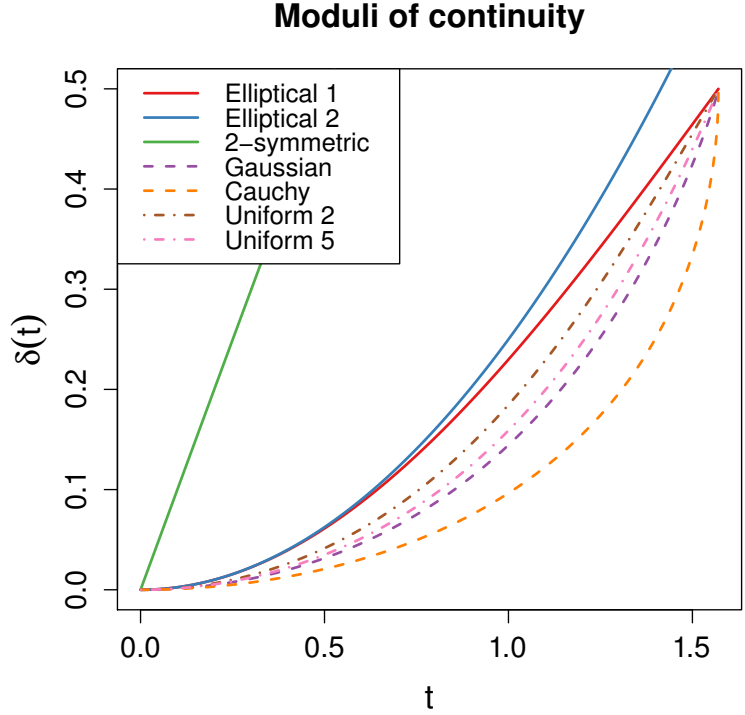}
	\caption{The moduli of continuity $\delta(t)$ of the halfspace functions $\varphi_x$ evaluated near their minima for various models. Elliptical 1 and 2 stand for the moduli from Theorem~\ref{theorem:elliptical} with $\delta_1$ and $\delta_2$, respectively; 2-symmetric is the modulus $\delta_1$ from Theorem~\ref{theorem:stable} for $p=2$; uniform 2 and 5 are the tight moduli for the uniform distribution with $d=2$ and $5$, respectively. The faster this function is decreases as $t \to 0$, the better approximation rate can be achieved.}
	\label{figure:delta}
\end{figure}

\begin{example}[Multivariate Gaussian distribution]
For a standard multivariate Gaussian distribution $X$, $\left\langle u, X \right\rangle$ is standard univariate Gaussian for any $u\in\Sph$, i.e.
	\[	F_p(t) = \int_{-\infty}^t \frac{1}{\sqrt{2 \pi}} \exp(-s^2/2) \dd s \quad\mbox{for }t\in\R.	\]
To get tight convergence rates, we must evaluate the modulus \eqref{delta} for $\varepsilon > 0$ fixed. It is not difficult to find that the maximal value of the latter function is attained at
	\[	t^* = \frac{\sqrt{-2 \log\left(\cos(\varepsilon)\right)}}{\sin(\varepsilon)} \quad \mbox{for }\varepsilon \in (0,\pi/2),	\]
which means that for $\left\Vert \widetilde{u}(x) - u \right\Vert_g \leq \varepsilon$ with $\varepsilon \in (0,\pi/2)$
	\begin{equation}	\label{tight bound}
	\sup_{x \in \R^d} \left\vert \varphi_x(\widetilde{u}(x)) - \varphi_x(u) \right\vert \leq \delta(\varepsilon) = F_p(t^*) - F_p(t^* \cos(\varepsilon)) \leq (1-\cos(\varepsilon))/2 \leq \varepsilon^2/4.	
	\end{equation}
The first bound above is the optimal one from Theorem~\ref{theorem:optimal}; the other two follow from Theorem~\ref{theorem:elliptical}. All of them can be used for numerical evaluation of the exact convergence rates. For that, it is enough to invert the inequality in Theorem~\ref{theorem:modulus} and obtain an error bound of the form $\Delta_n(C) \leq \delta \left( a_d^{-1} \left( \left( d \log \log n + \log n \right) / n \right) \right)$ which will approximately hold true almost surely for $n$ large enough. A comparison of the optimal uniform approximation errors for Gaussian distributions can be found in Table~\ref{table:all}. Graphically, these errors are compared with the general bounds from Theorem~\ref{theorem:elliptical} in Figure~\ref{figure:Gaussian_Cauchy}, see also the function $\delta(\varepsilon)$ from \eqref{tight bound}
displayed in Figure~\ref{figure:delta}. The tight bound presents a substantial improvement over the general bounds from Theorem~\ref{theorem:elliptical}.
\end{example}

\begin{table}[htpb]
\centering
\caption{Uniform error bounds for the approximation of the halfspace depth for elliptically symmetric unimodal distributions. Ellipt. 1 and 2 stand for the general bounds applicable to unimodal elliptically symmetric distributions from Theorem~\ref{theorem:elliptical} and $\delta_1$, $\delta_2$, respectively. 2-sym. is the general bound from Theorem~\ref{theorem:stable} for $p=2$ and $\delta_1$.} 
\label{table:all}
\begin{tabular}{cc|ccccc}
  \hline
 & $n$ & $d=2$ & $d=3$ & $d=5$ & $d=10$ & $d=20$ \\ 
  \hline
\parbox[t]{2mm}{\multirow{4}{*}{\rotatebox[origin=c]{90}{\footnotesize{Ellipt. 1}}}} & $10^2$ & 0.01441 & 0.09187 & 0.21855 & 0.35805 & 0.45576 \\ 
   & $10^3$ & 0.00029 & 0.01271 & 0.07629 & 0.20148 & 0.31099 \\ 
   & $10^4$ & 0.00000 & 0.00159 & 0.02625 & 0.11841 & 0.22779 \\ 
   & $10^5$ & 0.00000 & 0.00019 & 0.00892 & 0.07080 & 0.17139 \\  \hline 
\parbox[t]{2mm}{\multirow{4}{*}{\rotatebox[origin=c]{90}{\footnotesize{Ellipt. 2}}}}   & $10^2$ & 0.01448 & 0.09483 & 0.23663 & 0.41149 & 0.54924 \\ 
   & $10^3$ & 0.00029 & 0.01276 & 0.07831 & 0.21669 & 0.34995 \\ 
   & $10^4$ & 0.00000 & 0.00159 & 0.02648 & 0.12340 & 0.24756 \\ 
   & $10^5$ & 0.00000 & 0.00019 & 0.00895 & 0.07254 & 0.18219 \\  \hline 
\parbox[t]{2mm}{\multirow{4}{*}{\rotatebox[origin=c]{90}{\footnotesize{2-sym.}}}}   & $10^2$ & 0.24005 & 0.60619 & 0.93498 & 1.19675 & 1.35021 \\ 
   & $10^3$ & 0.03384 & 0.22544 & 0.55240 & 0.89774 & 1.11532 \\ 
   & $10^4$ & 0.00429 & 0.07968 & 0.32404 & 0.68821 & 0.95456 \\ 
   & $10^5$ & 0.00052 & 0.02745 & 0.18890 & 0.53218 & 0.82798 \\  \hline  \hline 
\parbox[t]{2mm}{\multirow{4}{*}{\rotatebox[origin=c]{90}{\footnotesize{Gaussian}}}}   & $10^2$ & 0.00707 & 0.04896 & 0.13535 & 0.27027 & 0.40894 \\ 
   & $10^3$ & 0.00014 & 0.00623 & 0.03997 & 0.12211 & 0.21854 \\ 
   & $10^4$ & 0.00000 & 0.00077 & 0.01305 & 0.06500 & 0.14277 \\ 
   & $10^5$ & 0.00000 & 0.00009 & 0.00436 & 0.03688 & 0.10010 \\  \hline 
\parbox[t]{2mm}{\multirow{4}{*}{\rotatebox[origin=c]{90}{\footnotesize{Cauchy}}}}   & $10^2$ & 0.00465 & 0.03226 & 0.09022 & 0.18834 & 0.31595 \\ 
   & $10^3$ & 0.00009 & 0.00410 & 0.02632 & 0.08119 & 0.14906 \\ 
   & $10^4$ & 0.00000 & 0.00051 & 0.00858 & 0.04288 & 0.09532 \\ 
   & $10^5$ & 0.00000 & 0.00006 & 0.00287 & 0.02428 & 0.06632 \\  \hline 
\parbox[t]{2mm}{\multirow{4}{*}{\rotatebox[origin=c]{90}{\footnotesize{Uniform}}}}   & $10^2$ & 0.00930 & 0.05831 & 0.14912 & --- & --- \\ 
   & $10^3$ & 0.00018 & 0.00743 & 0.04430 & --- & --- \\ 
   & $10^4$ & 0.00000 & 0.00092 & 0.01447 & --- & --- \\ 
   & $10^5$ & 0.00000 & 0.00011 & 0.00483 & --- & --- \\ 
   \hline
\end{tabular}
\end{table}

\begin{example}[Multivariate \texorpdfstring{$t$}{t}-distribution]
Consider the multivariate $t$-distribution $P\in\Prob$ with $\nu$ degrees of freedom, whose density is given by
	\[	f(x) = \frac{\Gamma\left(\frac{\nu+d}{2}\right)}{\Gamma(\nu/2) (\nu \pi)^{d/2} \left\vert \Sigma \right\vert^{1/2}}\left( 1 + \frac{1}{\nu}\left(x - \mu\right)\tr \Sigma^{-1} \left(x - \mu\right) \right)^{-(\nu+d)/2}	\quad\mbox{for }x\in\R^d,	\]
where $\Sigma \in \R^{d \times d}$ is non-singular, $\mu \in \R^d$, and $\nu=1,2,\dots$. For $\nu=1$ it is also known as the multivariate Cauchy distribution. Obviously, $P$ satisfies the assumptions of Theorem~\ref{theorem:elliptical}. 

We again focus on the case $\mu = 0$, and $\Sigma$ the identity matrix. The multivariate $t$-distribution is then spherically symmetric about the origin, with
	\[	F_p(t) = \int_{-\infty}^t \frac{\Gamma\left(\frac{\nu+1}{2}\right)}{\Gamma(\nu/2)\sqrt{\nu \pi}} \left( 1 + \frac{s^2}{\nu} \right)^{-(\nu+1)/2} \dd s \quad\mbox{for }t\in\R	\]
the distribution function of the univariate standard $t$-distribution with $\nu$ degrees of freedom, see \cite[Theorem~1]{Lin1972}.

Let us evaluate the tight modulus of continuity $\delta$ from \eqref{delta} as in the previous example. For general $\nu$, the argument of the maximum $t^*$ of the function $t \mapsto \left( F_p(t) - F_p(t\cos(\varepsilon)) \right)$ is not possible to be expressed in a closed form. For some special values of $\nu$ it can be computed that
	\[	t^* = \begin{cases}
						\left(\cos(\varepsilon)\right)^{-1/2} & \mbox{for $\nu = 1$}, \\
						\sqrt{\frac{3\left(-\left(\sin(\varepsilon)\right)^2-\left(\cos(\varepsilon )\right)^{3/2}+\sqrt{\cos(\varepsilon)}\right)}{\sqrt{\cos (\varepsilon )} \left(\left(\cos(\varepsilon)\right)^3-1\right)}} & \mbox{for $\nu=3$}, 
						\end{cases}	\quad\mbox{for }\varepsilon\in(0,\pi/2).	\]
For other small values of $\nu$ the formulas for $t^*$ can still be obtained explicitly, yet their expressions are already rather complicated. With $t^*$ known, the tight bound on the approximation errors \eqref{delta} can be obtained as in \eqref{tight bound}. In Figure~\ref{figure:Gaussian_Cauchy} and Table~\ref{table:all} we see the numerically computed bounds on the approximation error in the present setup for several values of $n$ and $d$. The function $\delta$ defined in \eqref{delta} can be found in Figure~\ref{figure:delta}. It can be seen that the approximation for heavy-tailed distributions $P$ is in fact more precise than for light-tailed ones. Indeed, if the tails of $P$ are rather flat and decrease slowly, the halfspace function $\varphi_x(u)$ will never change drastically for small deviations from its argument of the minimum $\widetilde{u}(x)$.  
\end{example}

\begin{figure}[htpb]
	\includegraphics[width=.45\textwidth]{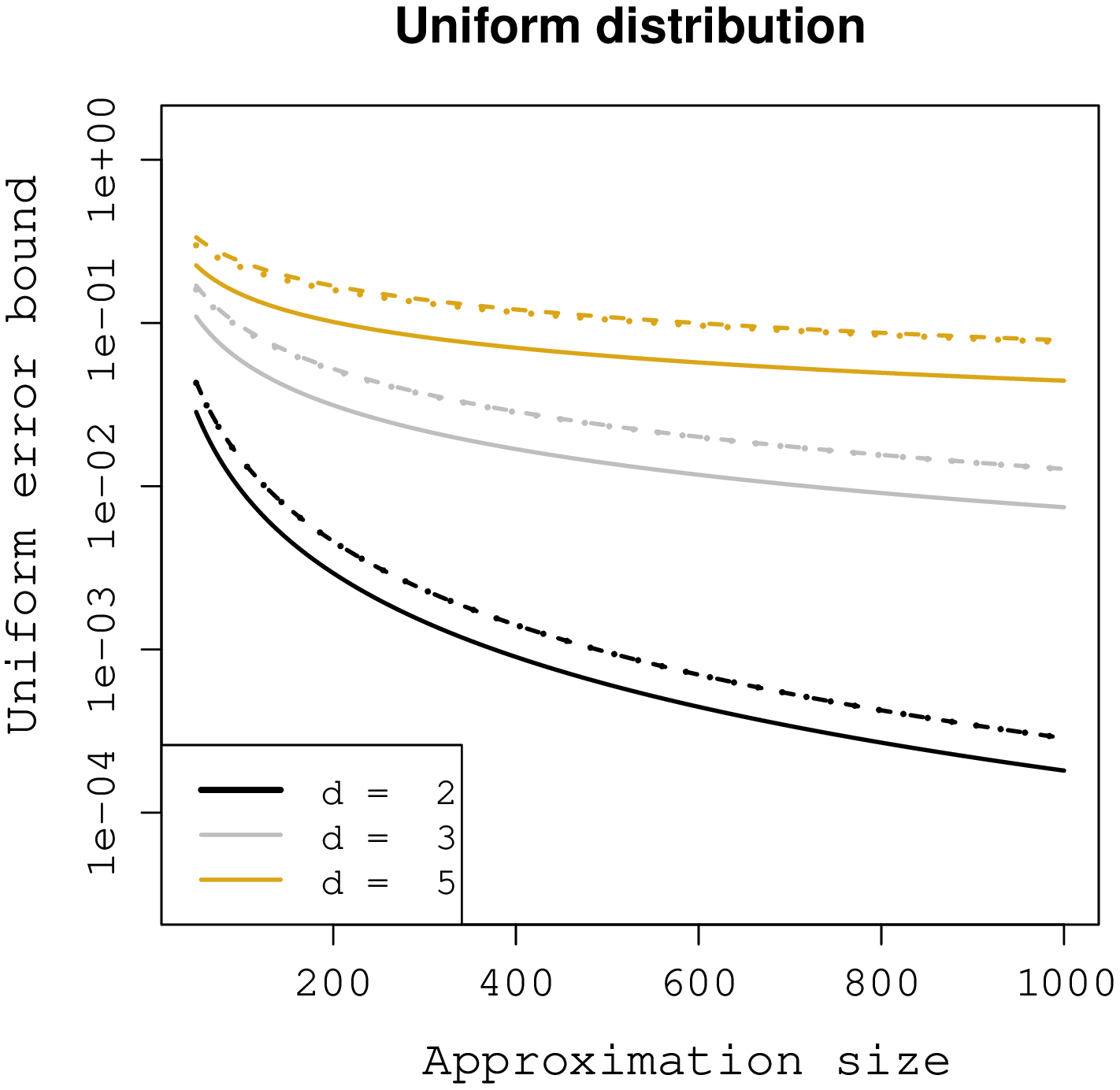} \qquad \includegraphics[width=.45\textwidth]{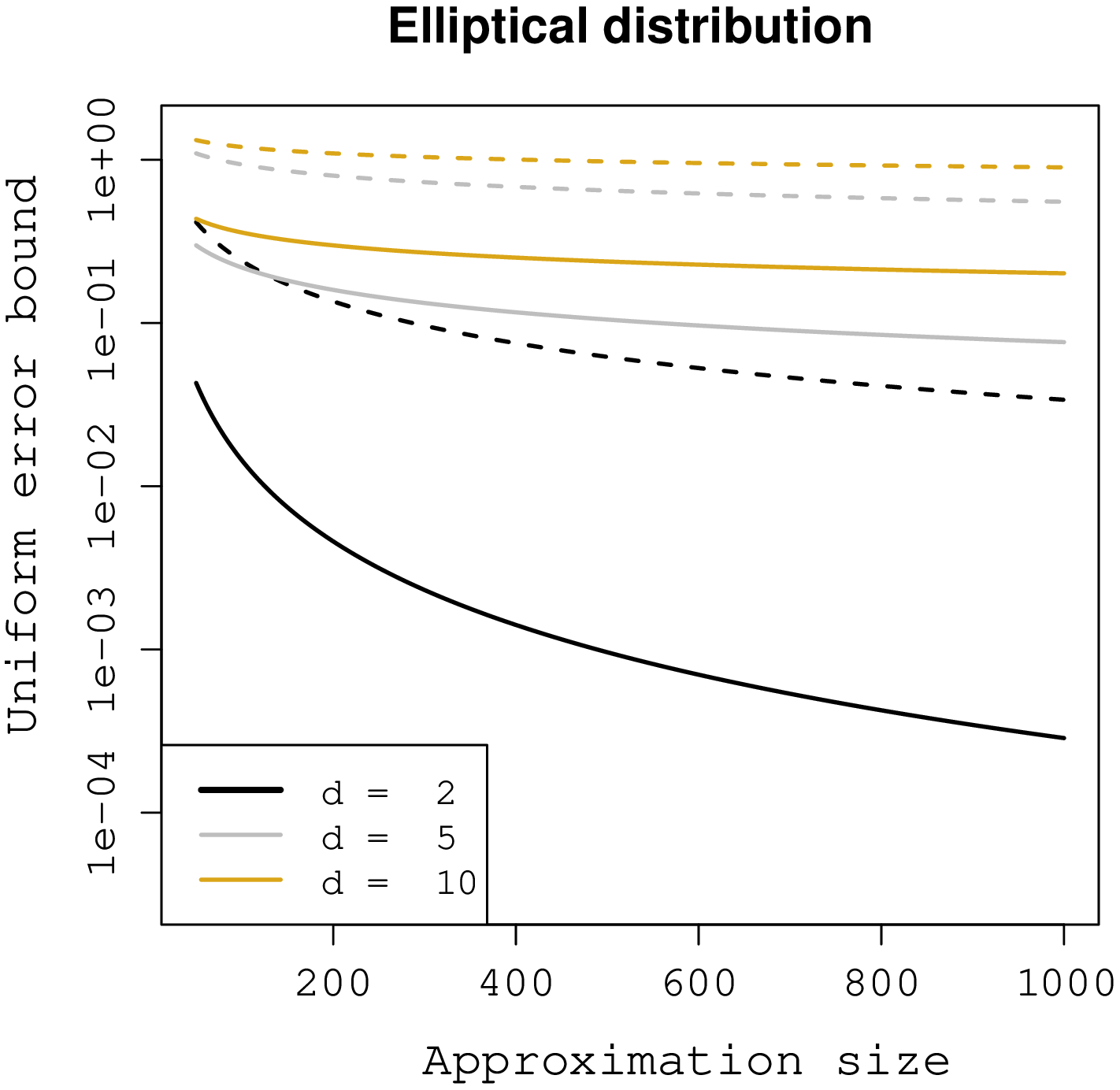}
	\caption{Uniform error bounds for $\Delta_n(\R^d)$ as a function of $n$ for the uniform distribution on the unit ball in $\R^d$ (left panel), and for general unimodal elliptically symmetric distributions (right panel). On the left panel, the solid curves represent the tight error bounds from Theorem~\ref{theorem:optimal}; the dashed and the dotted lines represent the two respective bounds from Theorem~\ref{theorem:elliptical}, valid for any unimodal elliptically symmetric distribution. The dashed and the dotted lines are nearly identical in the plot. On the right panel we compare bounds from Theorem~\ref{theorem:elliptical} with $\delta_1$ valid for unimodal elliptical distributions (solid lines), and the more general bounds from Theorem~\ref{theorem:stable} for $p=2$ and $\delta_1$ (dashed lines). Note the logarithmic scale of the vertical axes.}
	\label{figure:uniform_elliptical}
\end{figure}

\begin{example}[Uniform distribution on a ball]
Consider the uniform distribution $P$ on the unit ball in $\R^d$. The marginal distribution of $\left\langle u, X \right\rangle$ for any $u\in\Sph$ is given by the distribution function
	\[	F_p (t) = \int_{-1}^t \frac{\Gamma \left(d/2+1\right) \left(1-s^2\right)^{\frac{d-1}{2}}}{\sqrt{\pi} \, \Gamma \left(\frac{d+1}{2}\right)} \dd s \quad \mbox{for }t\in(-1,1).	\]
Note that unlike in our previous examples, $F_p$ depends on the dimension $d$. For small values of $d$, it is possible to express the point that realises the tight bound in \eqref{tight bound} explicitly
	\[	t^* = \begin{cases}
						\sqrt{\frac{2}{\cos (2 \varepsilon )+3}} & \mbox{for $d = 2$}, \\	
						\sqrt{\frac{\cos (\varepsilon )-1}{\left(\cos(\varepsilon )\right)^3-1}} & \mbox{for $d=3$}.
						\end{cases}	\]
In higher dimensions, approximations are easily obtained. Though, it is well known that as $d \to \infty$, the distribution of $\left\langle u, X \right\rangle$ is approximately Gaussian \citep{Diaconis_Freedman1984}, which reduces the problem for higher values of $d$ approximately to the case of the multivariate Gaussian distribution. In Table~\ref{table:all} and Figure~\ref{figure:uniform_elliptical} the error bounds are compared for $d=2$, $3$, and $5$. As can be seen, the approximation is slower than for both Gaussian and $t$-distributions. This is explained by the discontinuous nature of the density of $P$ --- for $x$ near the boundary of the unit ball and $d$ small, the halfspace function $\varphi_x$ deviates substantially from its minimum value already for $u$ rather close to its minimizer $\widetilde{u}(x)$. Though, the general bound from Theorem~\ref{theorem:elliptical} still performs somewhat worse.
\end{example}

\subsection{\texorpdfstring{$p$}{p}-symmetric distributions}	

A further important family of probability distributions whose halfspace depth can be evaluated explicitly are the $p$-symmetric distributions with $p \in (0,2]$, see \citep{Fang_etal1990}. That family of measures extends the spherically symmetric distributions treated in Section~\ref{section:elliptically symmetric distributions} that can be written as $p$-symmetric distributions with $p = 2$, as well as the multivariate extensions of stable distributions. Recall that for $p \in (0,2]$ a random vector $X = \left(X_1, \dots, X_d \right)\tr \sim P\in\Prob$ is said to have a $p$-symmetric distribution if for any $u \in \R^d$ the random variable $\left\langle X, u \right\rangle$ has the same distribution as $\left\Vert u \right\Vert_p X_1$, where $\left\Vert u \right\Vert_p = \left( \sum_{i=1}^d \left\vert u_i \right\vert^p \right)^{1/p}$ for $u = \left(u_1, \dots, u_d\right)\tr$, see \citep[Theorem~7.1]{Fang_etal1990}. We shall also write $\left\Vert u \right\Vert_\infty = \max_{i=1,\dots,d} \left\vert u_i \right\vert$. Observe that in this notation $\left\Vert \cdot \right\Vert_2 \equiv \left\Vert \cdot \right\Vert$.

For $P \in \Prob$ a $p$-symmetric distribution, the exact depth $\D(x)$ was computed in \cite[Example~(C)]{Masse_Theodorescu1994} and \cite[Theorem~3.1]{Chen_Tyler2004}. Let $F$ denote the (univariate) distribution function of the marginal $X_1$ of $X$, and set
	\begin{equation}	\label{q expression}
		  q = \begin{cases}
					\infty & \mbox{for }p\in(0,1], \\
					p/(p-1) & \mbox{for }p\in(1,2]
					\end{cases}	
	\end{equation}
to be the conjugate index to $p$. As demonstrated in \cite[Lemma~A.1 and Theorem~3.1]{Chen_Tyler2004}, we can write
	\begin{equation}	\label{depth for stable}
	\begin{aligned}
	\varphi_x(u) = F\left( \frac{\left\langle x, u \right\rangle}{\left\Vert u \right\Vert_p} \right), \quad \mbox{and} \quad	\D(x) = 1 - F\left(\left\Vert x \right\Vert_q\right) = F\left(-\left\Vert x \right\Vert_q\right).
	\end{aligned}
	\end{equation}
	
The following result generalizes Theorem~\ref{theorem:elliptical}. Its proof is in Appendix~\ref{appendix:stable}.
	
\begin{theorem}	\label{theorem:stable}
For $X = \left(X_1, \dots, X_d\right)\tr \sim P \in \Prob$ with a $p$-symmetric distribution such that the density of $X_1$ is unimodal we have that $\Delta_n(\R^d) \xrightarrow[n\to\infty]{\as}0$. Furthermore, the rate of convergence from Theorem~\ref{theorem:modulus} holds true with $C = \R^d$ and either
	\[
	\delta_1(t) = \begin{cases}
							d^{1/p-1/2} \left( d^{1/2-1/q} + 1\right) \sin\left(t/2\right) & \mbox{for }p\in[1,2],	\\
							d^{1/p}\sin\left(t/2\right) + \frac{d^{1/p-1/2}}{p} 2^{p-1}\sin^p\left(t/2\right) & \mbox{for }p\in(0,1),
							\end{cases}	\]
or the weaker
	\[
	\delta_2(t) = \begin{cases}
							d^{1/p-1/2} \left( d^{1/2-1/q} + 1\right) t/2 & \mbox{for }p\in[1,2],	\\
							d^{1/p} t/2 + \frac{d^{1/p-1/2}}{p} t^p/2 & \mbox{for }p\in(0,1),
							\end{cases}	\]
where in both expressions $t \in[0,\pi]$ and $q$ is given by \eqref{q expression}.
\end{theorem}

For $p = 2$ we can compare the obtained rate with that from Theorem~\ref{theorem:elliptical} and see that here the bounds are somewhat weaker. That is not surprising, since Theorem~\ref{theorem:stable} holds under much more general conditions. The unimodality assumption is satisfied for most commonly used $p$-symmetric distribution. For instance, it holds true for symmetric $p$-stable distributions \citep[Theorem~2.7.4]{Zolotarev1986}.

\begin{figure}[htpb]
	\includegraphics[width=.45\textwidth]{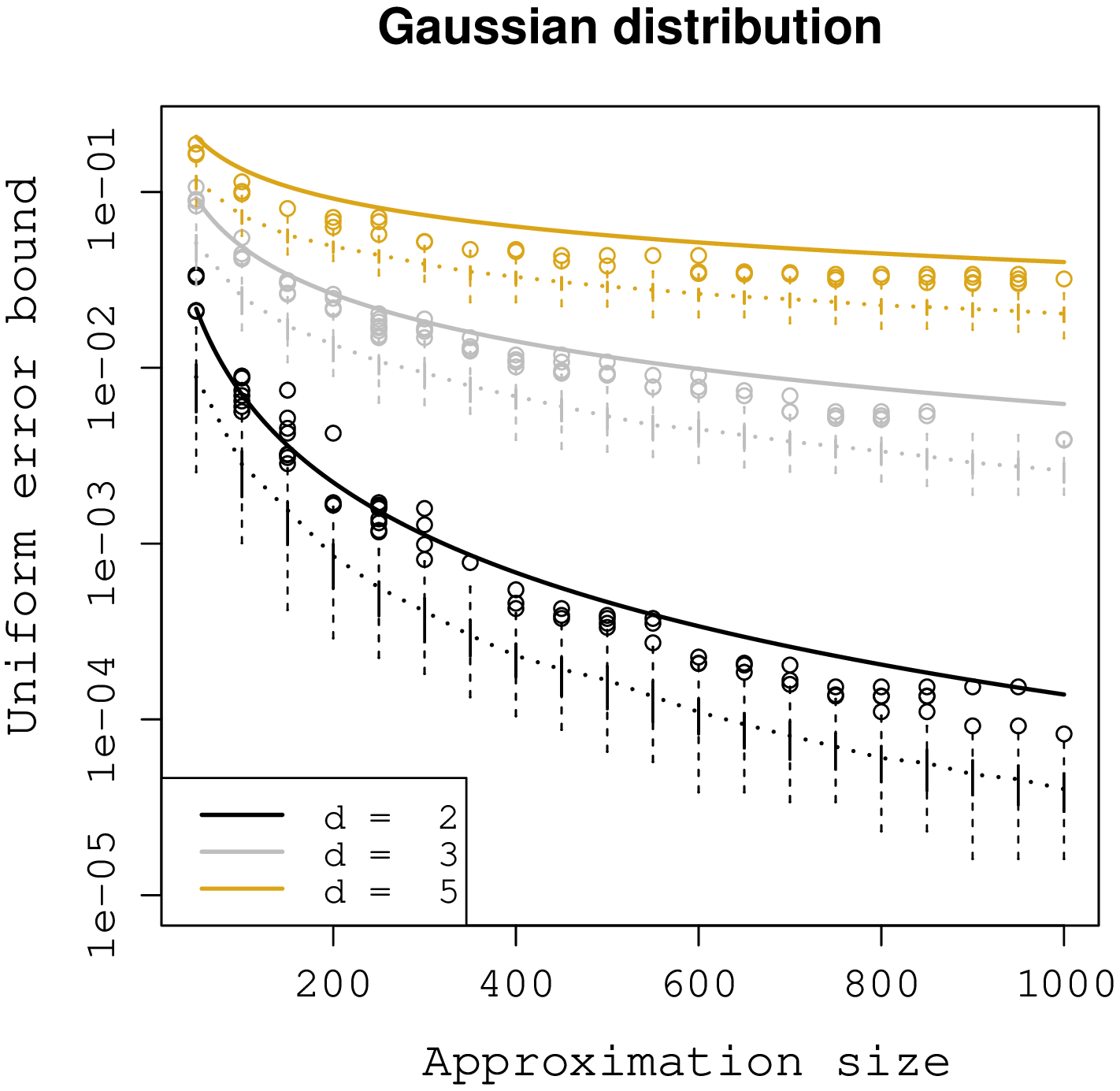} \qquad \includegraphics[width=.45\textwidth]{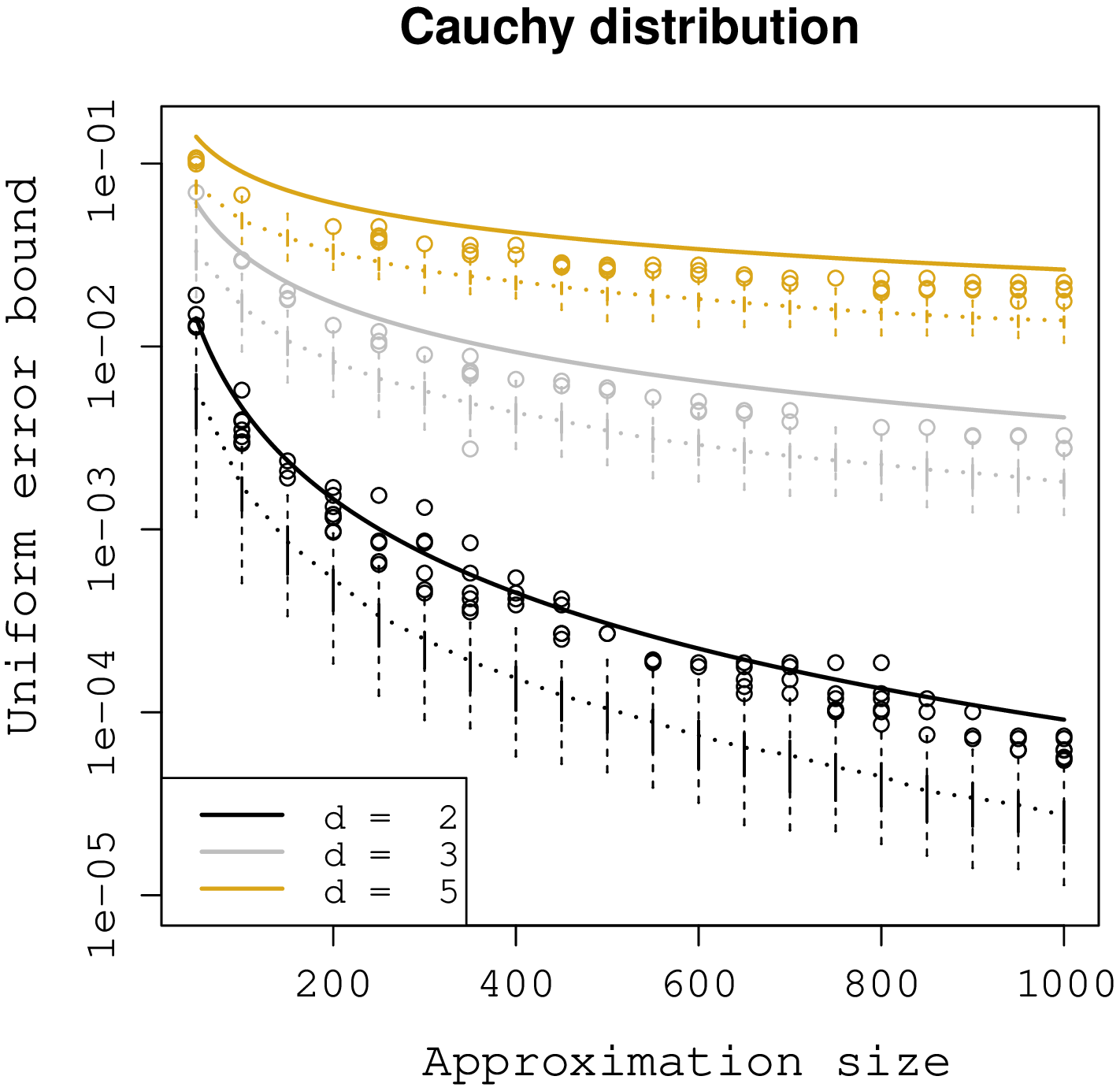} 
	\caption{The tight uniform bounds (solid lines) for $\Delta_n(\R^d)$ from Figure~\ref{figure:Gaussian_Cauchy} for the (standard) multivariate Gaussian distribution (left panel) and the (standard) multivariate Cauchy distribution (right panel). The bounds are complemented with boxplots of simulated data, and the dotted lines that interpolate the corresponding simulated means. In each of the $100$ independent runs, the supremum in $\Delta_n(\R^d)$ was approximated by the maximum value of $\left\vert \D(x_i\right) - \D_n\left(x_i) \right\vert$ for $500$ points $x_i$ sampled independently from the true distribution. The theoretical upper bounds match very well with the true approximation error.}
\label{figure:Gaussian_Cauchy_accuracy}
\end{figure}

\subsection{Non-uniform approximation}	\label{section:discontinuous case}

If the measure $P$ is not regular enough, uniform approximations of the halfspace depth may not be attainable. We illustrate our claim in two examples. In the first one we construct distributions that are atomic, yet their approximated depth fails to converge uniformly to its true counterpart. In the second example we demonstrate that if the density of $P$ is unbounded, uniform rates of convergence of the halfspace depth approximation are not possible to be derived. The latter result should be compared with Theorems~\ref{theorem:density example} and~\ref{theorem:bounded example}. 

\begin{example}
If $P$ does not admit a density, the halfspace functions $\varphi_x$ may contain discontinuities. In particular, this occurs when the depth is computed with respect to the empirical measure of a random sample of observations, i.e. for the sample halfspace depth. In that case, it can be inferred that the convergence of the approximations is not uniform, even on bounded sets $C$. 

Consider the example of $P$ atomic, supported in a set $\{x_i\}_{i=1}^m$ in $\R^d$ such that $P(\left\{x_i\right\}) = p_i$ with $\sum_{i=1}^m p_i = 1$. Let $0<p_1 \leq p_2 \leq \dots \leq p_m$. Without loss of generality, assume that the convex hull $K \subset \R^d$ of these points is $d$-dimensional, i.e. that $K$ is not fully contained in an affine subspace of $\R^d$ of dimension lower than $d$. Let $x$ be a point on a facet of $K$, and let $x_\varepsilon = x + \varepsilon \, v$ for $\varepsilon>0$, where $v \in \Sph$ is the outward normal of the facet of $K$ on which $x$ lies. Since $x_\varepsilon \notin K$ by definition, $\D(x_\varepsilon) = 0$ for any $\varepsilon>0$. For its approximation clearly $\D_n(x_\varepsilon) \geq p_1$ if and only if for some halfspace
	\[	H(x_\varepsilon,U_i) = \left\{ y \in \R^d \colon \left\langle U_i, y \right\rangle \leq \left\langle U_i, x_\varepsilon \right\rangle \right\}	\] 
whose boundary passes through $x_\varepsilon$ with outer normal $U_i \in \Sph$ intersects $K$. For $\varepsilon$ small, this will occur with high probability if the directions $U_i$ are sampled uniformly on $\Sph$ --- as $\varepsilon \to 0$ from the right, the condition $H(x_\varepsilon,U_i) \cap K = \emptyset$ effectively reduces to $U_i = -v$, an event of null probability. Thus, the convergence of the approximations cannot be uniform on the line segment $x + \varepsilon \, v$ with $\varepsilon \in (0,1)$, and
	\[	{\lim\sup}_{n\to\infty} \sup_{\varepsilon \in (0,1)} \left(\D_n(x_\varepsilon) - \D(x_\varepsilon)\right) \geq p_1 > 0 \quad\as \]
\end{example}

\begin{example}
Assume that $P$ is supported in a bounded subset of $\R^d$ and has a density, but that density is unbounded. By Theorem~\ref{theorem:density example} the uniform convergence of the approximations holds. But, no uniformly valid rate can be derived. To see this, we construct a distribution with an arbitrarily slow convergence rate. 

Let $\left\{x_i\right\}_{i=0}^\infty$ be distinct points on the upper half of the unit circle in $\R^2$, indexed in a clock-wise sense, and denote $l_i = \left\Vert x_{i} - x_{i-1} \right\Vert/2$ for $i\geq 1$, and set $l_0 = l_1$. Then, $\{ l_i \}_{i=0}^\infty$ is a sequence of positive numbers converging to zero. For $i \geq 0$, let $B_i$ denote a ball centred at $x_i$ with radius $e^{-i} l_i$, and let $C_i$ be the intersection of $B_i$ with the convex hull of the points $\{x_i\}_{i=0}^\infty$. Each $C_i$ is a convex subset of the unit ball in $\R^2$, with $x_i \in C_i$, and the sets $\{ C_i \}_{i=0}^\infty$ are pairwise disjoint. Let $P$ be the mixture of uniform distributions on the sets $C_i$ with mixing proportions $\{p_i\}_{i=0}^\infty$, where $\sum_{i=0}^\infty p_i = 1$ and $p_0 > p_1 > \dots > 0$. Consider the sequence of mid-points $y_j = (x_{j} + x_{j-1})/2$, $j=1,2,\dots$. Because each $y_j$ lies on the boundary of the convex hull of the support of $P$, the depth $\D(y_j)$ with respect to $P$ is zero for all $j$. Though, obviously, the single minimizer $\widetilde{u}(y_j)$ of the halfspace function $\varphi_{y_j}$ is the inwards normal of the facet corresponding to $x_{j-1}$ and $x_{j}$. If $\left\Vert u - \widetilde{u}(y_j) \right\Vert_g \geq \varepsilon$, then by the construction of $P$ for all $j \geq j(\varepsilon) = \left\lceil -\log \sin \varepsilon + 1 \right\rceil$, the smallest integer larger than $-\log \sin \varepsilon + 1$, it can be seen that either $C_{j}$, or $C_{j-1}$ lies inside the halfspace whose normal is $u$ and passes through $y_j$. Thus,
	\[	\varphi_{y_j}(u) = \PP\left(\left\langle u, X \right\rangle \leq \left\langle u, y_j \right\rangle \right) \geq p_{j}. \]
Using a result on maximal spacings in $\R$ from \cite[Theorem~5.2]{Devroye1981} it can be seen that for $U_1, \dots, U_n$ uniformly distributed on the circumference of the unit circle, it almost surely holds true that for all $n$ large enough and any $j$ we have $\min_{i=1,\dots,n} \left\Vert \widetilde{u}(y_j) - U_i \right\Vert_g \geq a(n)$, for a fixed positive sequence $\{a(n)\}_{n=1}^\infty$ that converges to zero. This means that almost surely for all $n$ large enough, for $j = \left\lceil - \log \sin a(n) + 1 \right\rceil$ we have
	\[	\sup_{x\in\R^2} \left\vert \D_n(x) - \D(x) \right\vert \geq \D_n(y_j) \geq p_j.	\]
Because $\{p_j\}_{j=0}^\infty$ can be made to converge to zero slowly, this means that no universal rate of convergence can be found, if the density of $P$ is allowed to be unbounded.
\end{example}

\section{Extensions to generalized projection depths}	\label{section:projection depth}

Let us now focus on another example of a depth that satisfies the projection property, yet is difficult to compute exactly --- the (generalized) projection depth. To define the depth, consider first the mappings
	\[
	\begin{aligned}
	m \colon \Prob[\R] \to \R \quad \mbox{and} \quad	s \colon \Prob[\R] \to [0,\infty],
	\end{aligned}
	\]
that satisfy the following conditions:
	\begin{itemize}
	\item $m(a X + b) = a\, m(X) + b$ for all $a, b \in \R$ and $X \sim P \in \Prob[\R]$;
	\item $s(a X + b) = a\, s(X)$ for all $a > 0$ and $X \sim P \in \Prob[\R]$. 
	\end{itemize}
Functionals $m$ and $s$ are called the \emph{location} and the \emph{scale parameter} of $X$, respectively. Using $m$ and $s$ for univariate distributions, it is possible to follow the ideas from \cite{Stahel1981, Donoho1982}, and define the outlyingness function in a multivariate space. For $x \in \R^d$ and $X \sim P \in \Prob$, the \emph{(projection) outlyingness} of $x$ with respect to $P$ is given by
	\begin{equation}	\label{outlyingness}
	O(x;P) = \sup_{u \in \Sph} \frac{\left\langle u, x \right\rangle - m(\left\langle u, X \right\rangle)}{s(\left\langle u, X \right\rangle)}.	
	\end{equation}
Function $O(\cdot;P)$ measures the largest deviation of a projection of a point from the location parameter of the corresponding projected distribution.

The outlyingness function \eqref{outlyingness} is closely related to depth --- high outlyingness indicates low centrality. Therefore, to construct a depth it suffices to transform the outlyingness index \eqref{outlyingness} by a non-increasing function $c$. A family of depth functions based on this idea is called the family of projection depths. It was studied in detail in \cite{Zuo_Serfling2000,Zuo2003}, and \cite{Dyckerhoff2004}. 

Consider a continuous function $c \colon [-\infty, \infty] \to [0,1]$ such that $c(x) = 1$ for $x<0$, and the restriction $c$ to $[0,\infty]$ is bijective, and strictly decreasing. The \emph{generalized projection depth} of $x \in \R^d$ with respect to $X \sim P \in \Prob[\R^d]$ is defined by
	\begin{equation}	\label{projection depth}
	\PD(x;P) = \inf_{u \in \Sph} c\left( \frac{\left\langle u, x \right\rangle - m(\left\langle u, X \right\rangle)}{s(\left\langle u, X \right\rangle)} \right).	
	\end{equation}
The family of depths \eqref{projection depth} was proposed in \cite{Dyckerhoff2004} as a generalization of the class of projection depths from \cite{Zuo2003}. The original projection depths are obtained by considering a scale functional $s$ invariant for reflections
	\[	s(-X) = s(X) \quad\mbox{for all }X \sim P\in\Prob[\R].	\]
Just as the halfspace depth, also the projection depth is difficult to calculate exactly, especially in high dimensions \cite{Zuo_Lai2011, Liu_Zuo2014}. In the same way as for the halfspace depth in \eqref{approximated depth}, define the approximated projection depth of $x \in \R^d$ with respect to $P \in \Prob$ based on the independent randomly sampled directions $U_1, \dots, U_n \in \Sph$ distributed uniformly on $\Sph$ by
	\[	\PD_n(x;P) = \PD_n(x) = \min_{i=1,\dots,n} c\left( \frac{\left\langle U_i, x \right\rangle - m(\left\langle U_i, X \right\rangle)}{s(\left\langle U_i, X \right\rangle)} \right).		\]
For a given set $C \subset \R^d$, we are interested in the uniform approximation of $\PD$ by $\PD_n$
	\[	\Delta_n^P(C) = \sup_{x \in C} \left\vert \PD_n(x) - \PD(x) \right\vert.	\]
Note that for the projection depth, the following function plays a role similar to that of the halfspace function from \eqref{phi function} for the halfspace depth
	\[	\varphi_x^P(u) = \frac{\left\langle u, x \right\rangle - m\left(\left\langle u, X \right\rangle \right)}{s \left(\left\langle u, X \right\rangle\right)}.	\]
A distinctive role will be now played by the argument(s) of the maxima of the function $\varphi_x^P(\cdot)$ on $\Sph$. Any representative of this set will be denoted by $\widetilde{u}^P(x)$. 
 
The next theorem provides simple conditions that guarantee the almost sure convergence of $\Delta^P_n(\R^d)$ to zero. In the proof of that result, the minimal multiplicative modulus of continuity of the function $c$ given by
	\begin{equation}	\label{zeta}
	\zeta(\tau) = \sup_{t > 0} \left( c(\tau \, t) - c(t) \right), \qquad\mbox{for }\tau \in (0,1),	
	\end{equation}
plays a crucial role. Note that for a continuous function $c$ the modulus $\zeta$ is a well-defined function continuous from the left at $\tau = 1$. To see this, pick a sequence $\tau_n \to 1$ from the left, and apply Dini's theorem \cite[Theorem~2.4.10]{Dudley2002} to the sequence of continuous functions $\psi_n(t) = c\left( \tau_n t \right)$ defined on a compact set $[0,\infty]$ (with $\psi_n(\infty) = 0$) that converge monotonically on $t \in [0,\infty]$ to $\psi_0(t) = c(t)$ (again, $\psi_0(\infty) = 0$). Dini's theorem asserts that this convergence must be uniform on the domain of $\psi_n$, which can be rewritten to $\zeta(\tau_n) = \sup_{t \geq 0} \left( \psi_n(t) - \psi_0(t) \right) = \sup_{t \geq 0} \left( c(\tau_n \, t) - c(t) \right) \xrightarrow[n\to\infty]{} 0$.

\begin{theorem}	\label{theorem:projection depth rate}
Let $X \sim P \in \Prob$ be such that the functions
	\[	
	\begin{aligned}
	u \mapsto m_u = m(\left\langle u, X \right\rangle), \quad \mbox{ and }\quad u \mapsto s_u = s(\left\langle u, X \right\rangle),
	\end{aligned}
	\]
are Lipschitz continuous on $\Sph$ with Lipschitz constants $C_m$ and $C_s$, respectively. Let $\sup_{u \in \Sph} \left\vert m_u \right\vert = \lambda_m < \infty$, $\inf_{u \in \Sph} s_u = \lambda_{s,1} > 0$, and $\sup_{u \in \Sph} s_u = \lambda_{s,2} < \infty$. Then $\Delta_n^P\left(\R^d\right) \xrightarrow[n\to\infty]{\as}0$. Furthermore, an explicit upper bound on the rate of convergence of $\Delta_n^P(\R^d)$ may be devised in terms of the constants above, the minimal modulus of continuity of $c$ given by
	\[	\delta_c(t) = \sup_{u,v > 0, \, \left\vert u - v \right\vert \leq t}	\left\vert c(u) - c(v) \right\vert,	\]
and the minimal multiplicative modulus of continuity of $c$ defined in \eqref{zeta}.
\end{theorem}

Theorem~\ref{theorem:projection depth rate} stands as an analogue of Theorem~\ref{theorem:Lipschitz} stated for the halfspace depth. Further minor modifications to the proof of Theorem~\ref{theorem:projection depth rate}, stated in Appendix~\ref{appendix:projection depth rate}, allow to formulate analogous results for continuous (thus uniformly continuous) location and scale parameters, in the spirit of Theorem~\ref{theorem:modulus}. We omit those technical details for brevity.

For the important case of spherically symmetric distributions, the rate from Theorem~\ref{theorem:projection depth rate} can be improved. The proof of this statement can be found in Appendix~\ref{appendix:projection depth spherical}.

\begin{theorem}	\label{theorem:projection depth spherical}
Let $P \in \Prob$ be a spherically symmetric distribution, and let $\zeta$ be a strictly increasing function whose inverse is denoted by $\zeta^{-1}$. Then
	\[
	{\lim\sup}_{n\to\infty} \frac{n \, a_d \left(1 - \zeta^{-1}\left(\Delta_n^P(\R^d)\right)\right) - \log n}{\log \log n} \leq d \quad\as	
	\]
Furthermore, for $c(x) = \frac{1}{1+(x)_+^k}$ with $k = 1, 2, \dots$ it holds true that
	\begin{equation}	\label{zeta function}
	\zeta(\tau) = \frac{1-\tau^{k/2}}{1+\tau^{k/2}} \qquad\mbox{for }\tau \in (0,1).	
	\end{equation}
\end{theorem}

Because the projection depths are affine invariant, analogously as in Theorem~\ref{theorem:affine} this result extends to elliptically symmetric distributions.

\subsection{Approximation of projection outlyingness}

Contrary to the projection depths \eqref{projection depth}, the outlyingness \eqref{outlyingness} is not uniformly approximated by its randomizations
	\[	O_n(x;P) = \max_{i=1,\dots,n} \frac{\left\langle U_i, x \right\rangle - m(\left\langle U_i, X \right\rangle)}{s(\left\langle U_i, X \right\rangle)},	\]
even if the conditions of Theorems~\ref{theorem:projection depth rate} and \ref{theorem:projection depth spherical} are satisfied. We illustrate this in a simple example. 

\begin{example}	
Let $X \sim P \in \Prob[\R^2]$ have the uniform distribution on the unit circle. Because each projection $\left\langle u, X \right\rangle$ has the same distribution symmetric about the origin, $m_u = 0$ and $s_u = S$ for all $u \in \SphO$ and some constant $S \geq 0$. For any reasonable scale functional $s$ the constant $S$ is positive. Consider a point $x = (x_1, 0)\tr \in \R^2$ with $x_1 > 0$. Then
	\[	\varphi_x^P(u) = \frac{\left\langle u, x \right\rangle - m_u}{s_u} = \frac{u_1 x_1}{S} \quad\mbox{for all }u = (u_1, u_2)\tr\in\SphO.	\]
This function is maximized at $\widetilde{u}^P(x) = (1, 0)\tr \in \Sph$ with $O(x;P) = x_1/S$. For any $u \in \SphO$ with $u_1 = \left\langle u, \widetilde{u}^P(x) \right\rangle < \cos(\varepsilon)$, on the other hand, $\varphi_x^P(\widetilde{u}^P(x)) - \varphi_x^P(u) > \left( 1 - \cos(\varepsilon) \right) x_1 /S$. Thus,
	\[	\lim_{n \to \infty} \sup_{x \in \R^d} \left\vert O(x;P) - O_n(x;P) \right\vert \geq \lim_{n \to \infty} \sup_{x_1 > 0} \max_{i = 1, \dots, n} \left(1 - \cos( \left\Vert \widetilde{u}^P(x) - U_i \right\Vert_g )\right) x_1 = \infty	\]
almost surely, because for any degree of approximation of the target direction $\widetilde{u}^P(x)$ there exists a point $x = \left(x_1, 0\right)\tr$ with a large norm so that its outlyingness is approximated poorly. 
\end{example}

\begin{figure}[htpb]
	\includegraphics[width=.65\textwidth]{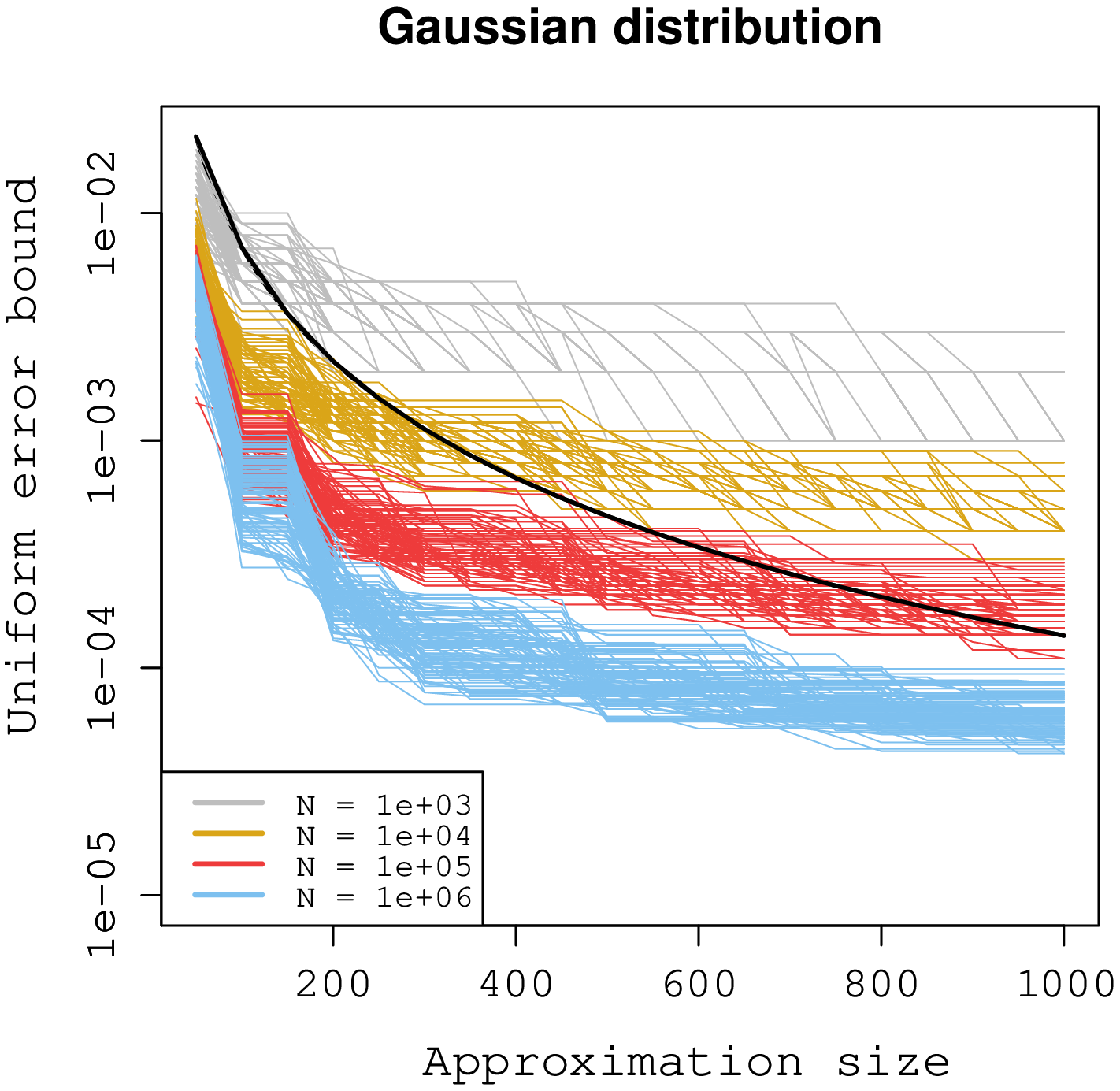}
	\caption{For a bivariate Gaussian distribution $P$, the tight bound from Theorem~\ref{theorem:optimal} on $\Delta_n(\R^2)$ (black solid line) is compared with $100$ independent realisations of simulated trajectories of $\max_{i=1,\dots,100} \left\vert \D_n\left(x_i, \widehat{P}_N\right) - \D\left(x_i, \widehat{P}_N\right) \right\vert$ as a function of the number of directions $n$, for $\left\{ x_i \right\}_{i=1}^N$ a random sample from $P$. Considered sample sizes are $N = 10^3, 10^4, 10^5$, and $10^6$.} 
	\label{figure:simulation}
\end{figure}

\section{Discussion} \label{section:other sampling distributions}

In the present paper we discussed the uniformity aspects of the depth approximation task. We demonstrated that for a distribution $P$ that is regular enough, the approximated depth does converge uniformly to its exact counterpart. This result justifies the simple approximation procedure from a theoretical point of view, as the almost sure uniform convergence of the depth approximations carries over to the approximated depth-statistics such as the depth median (i.e. an argument of the maximum of the depth on $\R^d$), or the level sets of the depth. 

In addition, we have presented and compared several almost sure upper bounds on the uniform discrepancies between the true halfspace and projection depth, and their approximated counterparts. Depending on the degree of regularity assumed about $P$, guidelines for the choice of the number of approximating directions $n$ in order to achieve a desired precision can be devised from the theory presented. For the halfspace depth, in Table~\ref{table:all} we saw that the random approximation scheme is feasible in lower dimensions, especially for distributions whose densities are continuous and rather flat. In dimensions $d>5$, hundreds of thousands of random directions may not be enough for sufficiently close approximations, and more elaborate algorithms appear to be needed. 

In Figure~\ref{figure:Gaussian_Cauchy_accuracy} we saw that for particular distributions, the theoretical bounds match the simulated results already for intermediate $n$ closely. Nevertheless, in practice, one does not observe $P$ directly, but rather only an empirical measure $\widehat{P}_N$ of a random sample $X_1, \dots, X_N$ of size $N = 1,2,\dots$ generated from $P$, and computes the depth with respect to $\widehat{P}_N$ as a surrogate for $\D(x;P)$ or $\PD(x;P)$. For practical considerations, it would therefore be desirable to obtain also upper bounds on the difference between the true depth and its approximations with respect to the empirical measure $\widehat{P}_N$. This can be done. For the halfspace depth, for instance, one may still use the theory presented, and devise a simple bound valid almost surely, for any fixed $\varepsilon > 0$ and all $N$ and $n$ large enough, that takes the form
	\[	
	\begin{aligned}
	\sup_{x\in C} & \left\vert \D_n(x;\widehat{P}_N) - \D(x;P) \right\vert \\
	& \leq \sup_{x\in C} \left\vert \D_n(x;\widehat{P}_N) - \D_n(x;P) \right\vert + \sup_{x\in C} \left\vert \D_n(x;P) - \D(x;P) \right\vert \\
	& \leq \sup_{x\in C} \max_{i=1,\dots,n} \left\vert \frac{1}{N} \sum_{j=1}^N \mathbb I\left[ \left\langle U_i, X_j \right\rangle \leq \left\langle U_i, x \right\rangle \right] - \varphi_x(U_i) \right\vert + \Delta_n(C) \\
	& \leq \sup_{x\in \R^d} \sup_{u\in\Sph} \left\vert \frac{1}{N} \sum_{j=1}^N \mathbb I\left[ \left\langle u, X_j \right\rangle \leq \left\langle u, x \right\rangle \right] - \varphi_x(u) \right\vert  + \Delta_n(C) \\
	& \leq \sqrt{\frac{\log \log N}{2 N}} + \varepsilon + \Delta_n(C).
	\end{aligned}
	\]
We used the fact that the sample halfspace depth process can be bounded by an empirical process given by the collection of closed halfspaces in $\R^d$. For the latter process, the bound follows from the law of the iterated logarithm devised in \cite[Corollary~2.4]{Kuelbs_Dudley1980}, see also \cite{Yukich1990}. The third summand in the last expression above can be handled using the results provided in this paper. Thus, the approximation of the empirical halfspace depth can be, at least for the number of observations $N$ large enough, approached by the approximation results for the true sampling distribution $P$, and an upper bound on the deviation of $P$ from $\widehat{P}_N$. An empirical argument supporting this finding is presented in Figure~\ref{figure:simulation}, where for $d=2$ and the standard Gaussian distribution, the theoretical bound from Theorem~\ref{theorem:optimal} is compared with simulated trajectories of $\Delta_{n,N} = \max_{i=1,\dots,100} \left\vert \D_n\left(x_i, \widehat{P}_N\right) - \D\left(x_i, \widehat{P}_N\right) \right\vert$ for several choices of $n \in [50, 1000]$ and four sample sizes $N = 10^3, 10^4, 10^5$ and $10^6$. Note that for sample size $N$, the quantity $\Delta_{n,N}$ is a multiple of $1/N$, and therefore unless $\Delta_{n,N} = 0$ (which did not happen in our simulation study), the simulated trajectory never decreases below $1/N$ (i.e. the tick $1e-3$ for $N = 1\,000$ etc.). In addition, because the supremum $\Delta_n(\R^2)$ is compared with a maximum of only $100$ points, also for $N$ large the simulated trajectories do not follow the theoretical bound as closely as in Figure~\ref{figure:Gaussian_Cauchy_accuracy}. Nevertheless, the obtained rates of convergence do appear to couple, even though the true distribution is replaced by an empirical one. Therefore, the general guidelines for the choice of $n$ are relevant also in the situation when the sample depth is approximated. For instance, if the theoretical bound in Table~\ref{table:all} is already too high for the practical application in mind, the simple approximation of the halfspace depth is certainly not a good idea, and more sophisticated methodologies must be employed. In any case, for the sample depth it must be kept in mind that according to the first example of Section~\ref{section:discontinuous case}, for empirical distributions the depth approximations are inherently non-uniform.

\appendix

\section{Proofs of the theoretical results: Halfspace depth}	\label{appendix:halfspace depth}

\subsection{Proof of Theorem~\ref{theorem:Lipschitz}}	\label{appendix:Lipschitz}
Denote by
	\begin{equation*} 
	S_n = \sup_{u\in\Sph} \min_{i=1,\dots,n} \left\Vert u - U_i \right\Vert_g
	\end{equation*}
the largest distance of a point on the sphere to the closest sample point $U_1,\dots,U_n$.
The quantity $S_n$ closely relates to the so-called maximal spacing problem, extensively studied in the literature. By \citet[the main theorem and the second remark on page~276]{Janson1987} and \cite[Theorem~1.2]{Janson1986} we know that in our setting it holds true that
	\begin{equation*}
	{\lim\sup}_{n\to\infty} \frac{n \, a_d(S_n) - \log n}{\log \log n} = d \quad\as, 
	\end{equation*}
where $a_d$ is the function defined in \eqref{a_d function}. Take $\varepsilon > 0$ fixed, and for all $x\in C$ let $\widetilde{u}(x) \in \Sph$ be such that 
	\[	\varphi_x(\widetilde{u}(x)) - \varepsilon < \D(x).	\]
By $j(x)$ denote any index such that
	\[	\left\Vert \widetilde{u}(x) - U_{j(x)} \right\Vert_g = \min_{i=1,\dots,n} \left\Vert \widetilde{u}(x) - U_i \right\Vert_g\le S_n.	\]
For the halfspace depth \eqref{halfspace depth} and its approximation \eqref{approximated depth} this implies that by \eqref{Lipschitz continuity}
	\begin{equation}	\label{rate derivation}
	\begin{aligned}
	0 \leq \Delta_n(C) & = \sup_{x \in C} \left( \D_n(x) - \D(x) \right) \\
	& < \sup_{x\in C} \left( \min_{i=1,\dots,n} \varphi_x(U_i) - \varphi_x(\widetilde{u}(x)) \right) + \varepsilon	\\
	& \leq \sup_{x\in C} \left( \varphi_x(U_{j(x)}) - \varphi_x(\widetilde{u}(x)) \right) + \varepsilon \\
	& \leq L \sup_{x\in C} \min_{i = 1, \dots, n } \left\Vert \widetilde{u}(x) - U_i \right\Vert_g + \varepsilon \\
	& \leq L \sup_{u\in\Sph} \min_{i = 1, \dots, n } \left\Vert u - U_i \right\Vert_g + \varepsilon \\
	& = L \, S_n + \varepsilon,
	\end{aligned}
	\end{equation}
where the last inequality holds for any $\varepsilon>0$. Therefore, $S_n\ge \Delta_n(C)/L$ and we may conclude that
	\[	{\lim\sup}_{n\to\infty}  \frac{n \, a_d\left( \Delta_n(C) / L\right) - \log n}{\log\log n} \leq d \quad\as	\]

\subsection{Proof of Theorem~\ref{theorem:modulus}}	\label{appendix:modulus}

The proof follows analogously to that of Theorem~\ref{theorem:Lipschitz}. The only difference is that in \eqref{rate derivation} we use a bound defined by means of the modulus of continuity $\delta$, instead of the Lipschitz property \eqref{Lipschitz continuity} of the class \eqref{phi functions}. Since $\|U_{j(x)} - \widetilde{u}(x)\|_g\le S_n$ and because of \eqref{modulus bound}, we obtain
	\begin{equation}	\label{basic inequality for general delta}
	\begin{aligned}
	0 \leq \Delta_n(C) & \leq \sup_{x\in C} \left( \varphi_x(U_{j(x)}) - \varphi_x(\widetilde{u}(x)) \right) + \varepsilon 
        \leq \delta\left( S_n \right) + \varepsilon
	\end{aligned}
	\end{equation}
for any $\varepsilon>0$, and the conclusion follows.

\begin{remark}
Note that the key step in the proofs of Theorems~\ref{theorem:Lipschitz} and~\ref{theorem:modulus} is a conceptually simple application of a continuity argument to an asymptotic result on maximal spacings in $\Sph$. Analogously, using the equicontinuity of the halfspace functions (or their analogues for the projection depth) only, any other asymptotic expression for the maximal spacings could be translated to an appropriate formula for the depth approximations using very similar proof techniques. We have opted for an asymptotic result of \citet{Janson1986, Janson1987} that takes the form of the law of the iterated logarithm for maximal spacings. This leads to upper bounds on the convergence rates that are quite strong in their formulation, yet may be rather conservative for some practical applications.
\end{remark}

\begin{remark}
If the random sample of approximating directions $U_1, \dots, U_n$ is taken from a non-uniform distribution on $\Sph$, results analogous to those described here can be obtained using the recent study of \citet{Aaron_etal2017}, which complements the theory of \citet{Janson1987}. 
\end{remark}

\subsection{Proof of Theorem~\ref{theorem:affine}} \label{appendix:affine}

For the approximated depth of the affinely transformed measure we can write
	\[	
	\begin{aligned}
	\D_n\left(A x + b; P_{A X + b}\right) & = \min_{i=1,\dots,n} \PP\left( \left\langle U_i, A X + b \right\rangle \leq \left\langle U_i, A x + b\right\rangle \right) \\
	& = \min_{i=1,\dots,n} \PP\left( \left\langle A\tr U_i, X \right\rangle \leq \left\langle A\tr U_i, x \right\rangle \right) \\
	&= \min_{i=1,\dots,n}\varphi_x\left(\frac{A\tr U_i}{\left\Vert A\tr U_i \right\Vert}\right)
	\end{aligned}
	\]
where on the right hand side we see the approximated depth of $x$ with respect to the untransformed random vector $X$, with the random directions $U_i$ sampled from a possibly non-uniform distribution on $\Sph$. As in \eqref{rate derivation} and \eqref{basic inequality for general delta} it follows that for any $\varepsilon > 0$
	\begin{align*}
        0\le\widetilde{\Delta}_n(C) 
        &= \sup_{x \in A C + b} \left\vert \D_n(x;P_{A X + b}) - \D(x;P_{A X + b}) \right\vert\\
        &= \sup_{x \in C} \left\vert \D_n(A x + b;P_{A X + b}) - \D(A x + b;P_{A X + b}) \right\vert\\
        & \leq \sup_{x \in C} \left\vert \min_{i=1,\dots,n}\varphi_x\left(\frac{A\tr U_i}{\left\Vert A\tr U_i \right\Vert}\right) - \varphi_x\left(\widetilde{u}(x)\right) \right\vert + \varepsilon.
        \end{align*} 
Thus, we need to bound the spherical maximal spacing that corresponds to a random sample from $A\tr U/\left\Vert A\tr U\right\Vert$ for $U$ uniform on $\Sph$. Let $u, v \in \Sph$. We bound the distance between the transforms of $u$ and $v$
	\begin{equation}	\label{affine spacing}
	\begin{aligned}
	\left\Vert \frac{A\tr u }{\left\Vert A\tr u \right\Vert} - \frac{A\tr v }{\left\Vert A\tr v  \right\Vert} \right\Vert_g & \leq \frac{\pi}{2}\left\Vert \frac{A\tr u }{\left\Vert A\tr u \right\Vert} - \frac{A\tr v }{\left\Vert A\tr v  \right\Vert} \right\Vert \\
	& \leq \frac{\pi}{2} \left( \left\Vert \frac{A\tr u }{\left\Vert A\tr u \right\Vert} - \frac{A\tr v }{\left\Vert A\tr u  \right\Vert} \right\Vert + \left\Vert \frac{A\tr v }{\left\Vert A\tr u \right\Vert} - \frac{A\tr v }{\left\Vert A\tr v  \right\Vert} \right\Vert \right) \\
	& = \frac{\pi}{2} \left( \frac{1}{\left\Vert A\tr u \right\Vert} \left\Vert A\tr \left( u - v \right) \right\Vert + \left\Vert A\tr v  \right\Vert \left\vert \frac{1}{\left\Vert A\tr u \right\Vert} - \frac{1}{\left\Vert A\tr v  \right\Vert} \right\vert \right) \\
	& \leq \frac{\pi}{2} \left( \frac{\lambda_1 \left\Vert u - v \right\Vert}{\lambda_d} + \frac{\lambda_1 \left\Vert A\tr u -  A\tr v \right\Vert}{\lambda_d^2} \right) \\
	& \leq \frac{\pi}{2} \left( \frac{\lambda_1 \left\Vert u - v \right\Vert}{\lambda_d} \left( 1 + \frac{\lambda_1}{\lambda_d} \right) \right)
        \leq K \left\Vert u - v \right\Vert_g  
	\end{aligned}
	\end{equation}
where $\lambda_1$ and $\lambda_d$ are the largest and the smallest singular value, respectively, of the matrix $A\tr$ and $K = \frac{\pi}{2} \frac{\lambda_1}{\lambda_d} \left( 1 + \frac{\lambda_1}{\lambda_d} \right)$. The third inequality above is justified by \cite[Theorem~3.1.2]{Horn_Johnson1994}, formula $\inf_{u \in \Sph} \left\Vert A\tr u \right\Vert = \lambda_d$ that follows from the same theorem, and the mean value theorem applied to function $g(t) = 1/t$. The first and the last inequality in \eqref{affine spacing} hold since for $u, v \in \Sph$ we have $\left\Vert u - v \right\Vert_g \leq \frac{\pi}{2} \left\Vert u - v \right\Vert \leq \frac{\pi}{2} \left\Vert u - v \right\Vert_g$. We can now bound the spherical maximal spacing $\widetilde{S}_n$ of the non-uniform distributed directions $A\tr U_i/\left\Vert A\tr U_i \right\Vert$ by
	\begin{align*}
	\widetilde{S}_n
	&= \sup_{u\in\Sph} \min_{i=1,\dots,n} \left\Vert u - \frac{A\tr U_i}{\left\Vert A\tr U_i \right\Vert} \right\Vert_g
= \sup_{v\in\Sph} \min_{i=1,\dots,n} \left\Vert \frac{A\tr v}{\left\Vert A\tr v \right\Vert} - \frac{A\tr U_i}{\left\Vert A\tr U_i \right\Vert} \right\Vert_g\\
	&\le \sup_{v\in\Sph} \min_{i=1,\dots,n} K \left\Vert v - U_i\right\Vert_g= K S_n.
\end{align*}
Now it follows that 
	\[
	0 \leq \widetilde{\Delta}_n(C) 
	\leq \delta\left(\widetilde{S}_n\right) + \varepsilon
        \leq \delta\left(K S_n\right) + \varepsilon
        \]
for any $\varepsilon > 0$, and the proof can be concluded as that of Theorem~\ref{theorem:modulus}.
	
%

\subsection{Proof of Theorem~\ref{theorem:density example}} \label{appendix:density example}

With no loss of generality, assume that $C \subset \R^d$ is compact with $P(C)=1$. Define the function
	\[	\varphi \colon C \times \Sph \to [0,1] \colon (x, u)\tr \mapsto \varphi_x(u).	\]
Because $P$ satisfies \eqref{Delta}, $\varphi$ is continuous on its domain \citep[Proposition~4.5]{Masse2004}, which is a compact subset of $\R^{2d}$ \citep[Theorem~2.2.8]{Dudley2002}. Therefore, $\varphi$ must be uniformly continuous on its domain \citep[Corollary~2.4.6]{Dudley2002} and for $\delta_\varphi$ the modulus of continuity of $\varphi$ defined by
	\[	\delta_\varphi \colon [0,\infty) \to [0,\infty) \colon t \mapsto \sup_{\left\Vert x - y \right\Vert + \left\Vert u - v \right\Vert_g \leq t} \left\vert \varphi(x,u) - \varphi(y,v) \right\vert	\]
it holds true that
	\[	\sup_{x \in C} \sup_{\left\Vert u - v \right\Vert_g \leq t} \left\vert \varphi_x(u) - \varphi_x(v) \right\vert = \sup_{x \in C} \sup_{\left\Vert u - v \right\Vert_g \leq t} \left\vert \varphi(x,u) - \varphi(x,v) \right\vert \leq \delta_\varphi(t) \quad\mbox{for all $t\geq0$}.	\]
Thus, we can apply Theorem~\ref{theorem:modulus} with the modulus $\delta_\varphi$, and the approximated depth approaches the theoretical one uniformly in $x$.

\subsection{Proof of Theorem~\ref{theorem:bounded example}}	\label{appendix:bounded example}

For any $x \in C$, $u, v \in \Sph$ we can bound
	\[	\left\vert \varphi_x(u) - \varphi_x(v) \right\vert \leq 2 \int_{W(r,u,v)} M \dd t,  	\]
where on the right hand side $W(r,u,v)$ is a wedge of a ball in $\R^d$ of radius $r=\mathrm{diam}(C)$ and center $x$ (that therefore contains $C$), whose bounding hyperplanes have normals $u, v \in \Sph$. Since the volume of this wedge is 
	\[	\left\Vert u - v \right\Vert_g \frac{\pi^{d/2-1} r^d}{2 \Gamma\left(d/2+1\right)}	\] 
we obtain
	\[
	\begin{aligned}
	\sup_{x\in C} \left\vert \varphi_x(u) - \varphi_x(v) \right\vert & \leq \frac{M \pi^{d/2-1} r^d}{\Gamma\left(d/2+1\right)} \left\Vert u - v \right\Vert_g,
	\end{aligned}
	\]
and Theorem~\ref{theorem:Lipschitz} provides the result.

\subsection{Proof of Theorem~\ref{theorem:elliptical}}	\label{appendix:elliptical}

Due to the translation invariance of $\D$, we can without loss of generality assume that $\mu = 0$. We thus assume that $P$ admits a density
	\begin{equation}	\label{density}
	f(x) = g\left( \left\Vert x \right\Vert^2 \right) \quad\mbox{for }x\in\R^d,
	\end{equation}
with $g \colon [0,\infty) \to [0,\infty)$. 
For such $P$, it is known that the true depth takes the form $\D(x) = h\left( \left\Vert x \right\Vert^2 \right)$ for some non-increasing function $h \colon [0,\infty) \to [0,\infty)$ \citep[Lemma~3.1]{Liu_Singh1993}. All contours of the depth $\D$ are therefore centred spheres. Consequently, for any $x \in \R^d$, $x \ne 0$, the single halfspace that realises the depth $\D(x)$ at $x$ is the one whose inner normal is $x /\left\Vert x \right\Vert$. In other words, for each $x \in \R^d$
	\[	\D(x) = \begin{cases}
							\varphi_x(u) & \mbox{for any $u\in\Sph$ if $x = 0$}, \\
							\varphi_x(-x/\left\Vert x \right\Vert) & \mbox{if $x\ne0$}.
							\end{cases}	\]
For $x = 0$, $\varphi_x(u) = 1/2$ for each $u\in\Sph$, and for any approximation we obtain the exact depth, i.e. $\D_n(x) - \D(x) = 0$ almost surely for all $n$. To assess the quality of the approximation at $x \ne 0$, it is necessary to control only the values $\varphi_x(u)-\varphi_x(\widetilde{u}(x))$ for $u$ close to $\widetilde{u}(x) = -x/\left\Vert x \right\Vert$. This makes the problem easier than the general case described in Theorems~\ref{theorem:Lipschitz} and~\ref{theorem:modulus} when no information about the position of $\widetilde{u}(x)$ is available. 

Since $X \sim P$ is assumed to be spherically symmetric about $0$, the random vector $A X$ has the same distribution as $X$ for any orthogonal matrix $A \in \R^{d \times d}$ \citep{Serfling_symmetry}. Thus, any projection $\left\langle u, X \right\rangle$ for $u \in \Sph$ has the same univariate distribution. Denote the cumulative distribution function of this projected random vector by 
	\[	F_p(t) = \PP\left( \left\langle u, X \right\rangle \leq t \right).	\]
Since $P$ was assumed to admit a density, 
$F_p$ is continuous. The halfspace functions $\varphi_x$ can be expressed in terms of $F_p$ 
	\begin{equation}	\label{intermediate1}
	\varphi_x(u) = \PP\left( \left\langle u, X \right\rangle \leq \left\langle u, x \right\rangle \right) = F_p\left( \left\langle u, x \right\rangle \right) \quad\mbox{for all $x\in\R^d$ and $u\in\Sph$}.	
	\end{equation}
For the depth, we can write
	\begin{equation}	\label{intermediate2}
	\D(x) = F_p\left( \left\langle \widetilde{u}(x), x \right\rangle \right) = F_p\left(-\left\Vert x \right\Vert \right) \quad\mbox{for all $x\in\R^d$}.	
	\end{equation}
Writing again $U_{j(x)}$ for the closest element of the random sample $U_1, \dots, U_n$ from $\widetilde{u}(x)$, now we can bound
	\begin{equation}	\label{uniform limit}
	\begin{aligned}
	\sup_{x \in \R^d} \left\vert \D_n(x)-\D(x)\right\vert 
        & \le \sup_{x\in \R^d} \left\vert \varphi_x(U_{j(x)}) - \varphi_x(\widetilde{u}(x)) \right\vert \\ 
        & = \sup_{x\in\R^d} \left\vert F_p\left(\left\langle U_{j(x)}, x \right\rangle\right) - F_p\left(-\left\Vert x \right\Vert\right) \right\vert \\
	& = \sup_{x\in\R^d} \left\vert F_p\left(- \left\Vert x \right\Vert \left\langle U_{j(x)}, \widetilde{u}(x) \right\rangle\right) - F_p\left(-\left\Vert x \right\Vert\right) \right\vert \\
	& = \sup_{x\in\R^d} \left\vert F_p\left(- \left\Vert x \right\Vert \cos\left(\left\Vert U_{j(x)} - \widetilde{u}(x)\right\Vert_g\right)\right) - F_p\left(-\left\Vert x \right\Vert\right) \right\vert \\
	& \le \sup_{t \geq 0} \left\vert F_p\left( -t \cos\left(S_n\right) \right) - F_p(-t)\right\vert \\
	& = \sup_{t \geq 0} \left( F_p(t) - F_p(t\cos(S_n)) \right),	
	\end{aligned}
	\end{equation}
where we used the fact that $\left\Vert U_{j(x)} - \widetilde{u}(x)\right\Vert_g\le S_n$, further that $F_p$ is non-decreasing, and that the distribution of $\left\langle u, X \right\rangle$ is symmetric about zero, i.e. $F_p(t)=1-F_p(-t)$. Take an arbitrary decreasing sequence $\left\{ \varepsilon_\nu \right\}_{\nu=1}^\infty$ such that $\varepsilon_\nu \xrightarrow[\nu\to\infty]{}0$, and define a sequence of functions
	\begin{equation}	\label{F sequence}
		  \left\{ F_\nu(t) = \begin{cases}
												 F_p\left(t \cos\left(\varepsilon_\nu\right)\right) & \mbox{for }t \in [0,\infty), \\
												 1 & \mbox{for }t=\infty.
												\end{cases} \right\}_{\nu=1}^\infty.	
	\end{equation}
This sequence consists of continuous functions defined on a compact set $[0,\infty]$. As $F_p$ is continuous and non-decreasing,
	\[
	\begin{aligned}
	F_\nu(t) & \xrightarrow[\nu\to\infty]{} F_p(t) \\
	F_\nu(t) & \leq F_{\nu+1}(t)
	\end{aligned} \quad\mbox{for all $t\in[0,\infty]$ and $\nu = 1,2,\dots$}
	\]
Therefore, it is possible to use Dini's theorem \citep[Theorem~2.4.10]{Dudley2002} to obtain that the sequence of functions \eqref{F sequence} converges to $F_p$, uniformly in $[0,\infty]$. Thus, the limit of the right hand side of \eqref{uniform limit} is zero and finally, using a result on maximal spacings of \citet{Janson1987} we may conclude that the uniform convergence of the approximated depth holds true.
	
Let us now obtain the rate of convergence of the approximation if the density $f$ is unimodal. First, we show that for $d>1$ and $t\geq 0$ the distribution function $F_p(t)$ must be concave. To see this, note that for the density of the random variable $\left\langle u, X \right\rangle$ with $u = \left(0, \dots, 0, 1\right) \in \Sph$ we have
	\[	f_p(t) = \int_{\R} \dots \int_{\R} g \left(\left\Vert \left(x_1, \dots, x_{d-1}, t\right) \right\Vert^2 \right) \dd x_1 \dots \dd x_{d-1} \quad \mbox{for $t \in \R$}.	\]
For $f$ unimodal, function $g$ in \eqref{density} must be non-increasing. For $0 \leq t_1 \leq t_2$ this means that $g \left(\left\Vert \left(x_1, \dots, x_{d-1}, t_1\right) \right\Vert \right) \geq g \left(\left\Vert \left(x_1, \dots, x_{d-1}, t_2\right) \right\Vert \right)$ for any $x_1, \dots, x_{d-1} \in \R$, and that
	\[	f_p(t_1) \geq f_p(t_2) \quad\mbox{for any $0 \leq t_1 \leq t_2$}.	\]
Since $f_p(t)$ is non-increasing for $t\geq 0$, the distribution function $F_p(t)$ must be concave for $t \geq 0$. Thus, using the obvious fact that $F_p(0) = 1/2$, for $S_n\le \pi/2$ and $t \geq 0$ we have
	\begin{equation}	\label{concavity bound}
	\begin{aligned}
	0 & \leq F_p(t) - F_p(\cos(S_n)t) \leq F_p(t) - (1-\cos(S_n)) F_p(0) - \cos(S_n) F_p(t) \\
	& = \left(1-\cos(S_n)\right) F_p(t) - (1-\cos(S_n))/2 \leq (1-\cos(S_n))/2.	
	\end{aligned}
	\end{equation}
Combining the last formula with \eqref{uniform limit} gives 
	\begin{equation}	\label{cosine expression}
	\sup_{x \in \R^d} \left\vert \D(x) - \D_n(x) \right\vert \leq (1-\cos(S_n))/2.	
	\end{equation}
This inequality holds true under the condition that for any $u \in \Sph$ there exists a random direction $U_i$ such that $\left\Vert u - U_i \right\Vert_g \leq S_n$. In other words, the inequality remains valid if the polar angle that corresponds to the maximal spacing given by the random sample of directions $U_1, \dots, U_n$ does not exceed $S_n$. This leads to an inequality of the same type as \eqref{basic inequality for general delta} from Theorem~\ref{theorem:modulus}, and the conclusion follows by the same argument as that used in the proof of Theorem~\ref{theorem:modulus}. The simpler bound follows by application of the inequality $1 - \cos(t) \leq t^2/2$ to \eqref{cosine expression}. This way one obtains
	\[	\sup_{x \in \R^d} \left\vert \D(x) - \D_n(x) \right\vert \leq S_n^2/4,	\]
for $S_n$ the maximal spacing from the proofs of Theorems~\ref{theorem:Lipschitz} and~\ref{theorem:modulus}. Again, the same technique as that from the proof of Theorem~\ref{theorem:modulus} gives the final rate of convergence.

\subsection{Proof of Theorem~\ref{theorem:optimal}}	\label{appendix:optimal}

Again, without loss of generality $\mu = 0$. It is enough to show the first inequality --- the second assertion then follows from Theorem~\ref{theorem:elliptical}. Consider the main theorem of \citet{Janson1987} once more. By the first part of that result, for $S_n$ the spacings defined in the proof of Theorem~\ref{theorem:Lipschitz} we have
	\begin{equation}	\label{R expression}
	{\lim\inf}_{n\to\infty} \frac{n \, a_d(S_n) - \log n}{\log \log n} = d-2 \quad\as, 
	\end{equation}
which means that almost surely for $n$ large there exists a vector $v_n \in \Sph$ such that each sampled direction $U_i$, $i=1,\dots,n$, is far enough from $v_n$, i.e.  for any $\varepsilon > 0$ almost surely for all $n$ large enough there is $v_n$ such that
	\begin{equation}	\label{v_n inequality}
	\min_{i=1,\dots,n} \left\Vert v_n - U_i \right\Vert_g \geq a_d^{-1}\left(\frac{\left( d - 2 - \varepsilon \right) \log \log n + \log n}{n} \right) - \varepsilon.	
	\end{equation}
Denote the right hand side of the previous inequality by $r_n$.

By the definition of $\delta$ from \eqref{delta} for any $\eta>0$ there exists $K_n > 0$ such that
	\[	F_p(K_n) - F_p(K_n \cos(r_n) ) = F_p(-K_n \cos(r_n) ) - F_p(-K_n) \geq \delta(r_n) - \eta.	\]
Consider the point $x_n = -K_n v_n$. Then $\left\Vert x_n \right\Vert = K_n$, and $v_n$ is the minimizer of $\varphi_{x_n}$, i.e. $v_n = \widetilde{u}(x_n)$. The latter follows from the proof of Theorem~\ref{theorem:elliptical}, where it was argued that for a spherically symmetric $P$, $\widetilde{u}(x_n)$ minimizes $\varphi_{x_n}$ for all $x_n$ on the ray $\lambda \, \widetilde{u}(x_n)$ with $\lambda<0$. Let $u \in \Sph$ be any direction with $\left\Vert u - \widetilde{u}(x_n) \right\Vert_g \geq r_n$. By formulas \eqref{intermediate1} and \eqref{intermediate2} from the proof of Theorem~\ref{theorem:elliptical} and the monotonicity of $F_p$ we have
	\[
	\begin{aligned}
	\varphi_{x_n}(u) - \varphi_{x_n}\left(\widetilde{u}(x_n)\right) & = F_p\left(-\left\Vert x_n \right\Vert \cos\left(\left\Vert u - \widetilde{u}(x_n) \right\Vert_g\right)\right) - F_p(-\left\Vert x_n \right\Vert) \\
	& \geq F_p(-\left\Vert x_n \right\Vert \cos(r_n)) - F_p(-\left\Vert x_n \right\Vert) \geq \delta(r_n) - \eta.
	\end{aligned}
	\]
We have that if there exists $v_n \in \Sph$ that satisfies \eqref{v_n inequality}, then for any $\eta > 0$ it is possible to find $x_n$ such that for random directions $U_i$ the above inequality holds true with $u$ replaced by $U_i$ for all $i = 1, \dots, n$. In other words, we have that in this situation
	\begin{equation}	\label{delta formula}
	\Delta_n(\R^d) \geq \D_n(x_n) - \D(x_n) = \min_{i=1,\dots,n} \left( \varphi_{x_n}(U_i) - \D(x_n) \right) \geq \delta\left(r_n\right) - \eta.
	\end{equation}	
By \eqref{R expression} and \eqref{v_n inequality}, however, we know that almost surely for all $n$ large such $v_n$ exists. Therefore, it is possible to invert \eqref{delta formula}, which allows us to write that almost surely for any $\varepsilon > 0$ and $\eta>0$
	\[	{\lim\inf}_{n\to\infty} \frac{n \, a_d \left(\delta^{-1}\left(\Delta_n(\R^d) + \eta \right) + \varepsilon \right) - \log n}{\log \log n} \geq d - 2 - \varepsilon,	\]
and the conclusion follows.

\subsection{Proof of Theorem~\ref{theorem:stable}}	\label{appendix:stable}

The proof is led in a spirit similar to that of Theorem~\ref{theorem:elliptical}. Define
	\[	r \colon \R^d\setminus\{0\} \times \Sph \to \R \colon \left(x\tr, u\tr\right)\tr \mapsto \left. \frac{\left\langle x, u \right\rangle}{\left\Vert u \right\Vert_p}\middle/\frac{\left\langle x, \widetilde{u}(x) \right\rangle}{\left\Vert \widetilde{u}(x) \right\Vert_p}\right..	\]
For $p \leq 1$ the minimizing direction $\widetilde{u}(x) \in \Sph$ is not uniquely defined. Nevertheless, from \eqref{depth for stable} we see that for any choice of the minimal direction
	\begin{equation}	\label{q norm argument}
	\frac{\left\langle x, \widetilde{u}(x) \right\rangle}{\left\Vert \widetilde{u}(x) \right\Vert_p} = - \left\Vert x \right\Vert_q,	
	\end{equation}
and the function $r$ is well defined with
	\[	r(x,u) = -\frac{\left\langle x, u \right\rangle}{\left\Vert x \right\Vert_q \left\Vert u \right\Vert_p}.	\]
The essential part of this proof is to show that the function $r$ converges to $1$ uniformly in $x \in \R^d\setminus\{0\}$ when $u \in \Sph$ comes from a small neighbourhood of $\widetilde{u}(x)$. More precisely, we show that for any $\varepsilon>0$ there exists $t > 0$ such that
	\begin{equation}	\label{s uniform convergence}
	 \sup_{x \in \R^d \setminus\{0\}} \sup_{\substack{u \in \Sph \\ \left\Vert u - \widetilde{u}(x) \right\Vert_2 \leq t}} \left\vert r(x,u) - 1 \right\vert \leq \varepsilon.
	\end{equation}
Note that for the special case $p = 2$ it is possible to proceed as in the proof of Theorem~\ref{theorem:elliptical} and write (see formula~\eqref{uniform limit})
	\[	\left\vert r(x,u) - 1 \right\vert = \left\vert \frac{\left\langle x, u \right\rangle}{\left\langle x, \widetilde{u}(x) \right\rangle} - 1 \right\vert = \left\vert \left\langle u, \widetilde{u}(x) \right\rangle - 1 \right\vert = \left\vert 1 - \cos\left(\left\Vert u - \widetilde{u}(x) \right\Vert_g \right) \right\vert.	\]
For general $p \in (0,2]$ the derivation below is somewhat similar, yet more involved. 

From formula \eqref{s uniform convergence} it follows that for any $x \in \R^d \setminus \{ 0 \}$ and $u \in \Sph$ with $\left\Vert u - \widetilde{u}(x) \right\Vert_2 \leq t$ we can write
	\[	\frac{\left\langle x, u \right\rangle}{\left\Vert u \right\Vert_p} \in \left[ -(1+\varepsilon) \left\Vert x \right\Vert_q, - (1-\varepsilon) \left\Vert x \right\Vert_q \right], 	\]
and
	\[	\varphi_x(u) = F\left(\frac{\left\langle x, u \right\rangle}{\left\Vert u \right\Vert_p}\right) \in \left[ F\left( -(1+\varepsilon) \left\Vert x \right\Vert_q \right), F\left( - (1-\varepsilon) \left\Vert x \right\Vert_q \right)\right].	\]
By the unimodality assumption, the distribution function $F$ must be concave on $[0,\infty)$. Thus, using derivation analogous to \eqref{concavity bound} we can write
	\[	
	\begin{aligned}
	0 & \leq \varphi_x(u) - \varphi_x\left(\widetilde{u}(x)\right) = \varphi_x(u) - \D(x) = F\left( \frac{\left\langle x, u \right\rangle}{\left\Vert u \right\Vert_p} \right) - F\left(-\left\Vert x \right\Vert_q\right) \\
	& \leq F\left( - (1-\varepsilon) \left\Vert x \right\Vert_q \right) - F\left(-\left\Vert x \right\Vert_q\right) = F\left(\left\Vert x \right\Vert_q\right) - F\left( (1-\varepsilon) \left\Vert x \right\Vert_q \right) \leq \frac{\varepsilon}{2}.	
	\end{aligned}
	\]
The last inequality holds for all $x \ne 0$ and $u \in \Sph$ with 
	\begin{equation}	\label{g 2 norm relation}
	\left\Vert u - \widetilde{u}(x) \right\Vert_2 = 2\sin\left(\left\Vert u - \widetilde{u}(x) \right\Vert_g /2\right) \leq t.	
	\end{equation}
Therefore, using the same argumentation as in the proof of Theorem~\ref{theorem:modulus} we can conclude that the uniform convergence of the approximations holds true as claimed.

It remains to show \eqref{s uniform convergence}. First, observe that for any $\alpha>0$ it holds that $r(x,u) = r(\alpha \, x, u)$. Thus, in \eqref{s uniform convergence} we may focus only on $x \in \Sph$ in the first supremum. Assume that $u \in \Sph$ is such that $\left\Vert u - \widetilde{u}(x) \right\Vert_2 \leq t$. 
Let us rewrite $r$ as a product of two terms
	\begin{equation}	\label{r ratio}
	r(x,u) = \frac{\left\Vert \widetilde{u}(x) \right\Vert_p}{\left\Vert u \right\Vert_p} \frac{\left\langle x, u \right\rangle}{\left\langle x, \widetilde{u}(x) \right\rangle}	
	\end{equation}
and inspect them separately. For the first factor we have to distinguish the cases $p\in(0,1)$ and $p\in [1,2]$.

The following relations hold true for $0<p < q \leq \infty$ and $x \in \R^d$ \citep[see, e.g.,][Lemma~A.1]{Chen_Tyler2004}
	\begin{equation}	\label{l_p relations}
	\left\Vert x \right\Vert_q \leq \left\Vert x \right\Vert_p \leq d^{1/p-1/q} \left\Vert x \right\Vert_q.	
	\end{equation}
With $q=2$ and since $\|u\|_2=\|\widetilde{u}(x)\|_2=1$ it follows that
	\[
        \frac{\|\widetilde{u}(x)\|_p}{\|u\|_p}
        \le \frac{d^{1/p-1/2}\|\widetilde{u}(x)\|_2}{\|u\|_2}
        =d^{1/p-1/2}\,.
        \]	
For $p\geq 1$ the function $\left\Vert \cdot \right\Vert_p$ is a norm. Thus, it obeys the (reverse) triangle inequality, and for the first factor in \eqref{r ratio} we can bound for $p\in[1,2]$
	\begin{equation}	\label{term 1}
	\left\vert \frac{\left\Vert \widetilde{u}(x) \right\Vert_p}{\left\Vert u \right\Vert_p} - 1 \right\vert = \frac{1}{\left\Vert u \right\Vert_p} \left\vert \left\Vert \widetilde{u}(x) \right\Vert_p - \left\Vert u \right\Vert_p \right\vert \leq \frac{\left\Vert u - \widetilde{u}(x) \right\Vert_p}{\left\Vert u \right\Vert_p} \leq d^{1/p-1/2} \left\Vert u - \widetilde{u}(x) \right\Vert_2.
	\end{equation}
Now, we discuss the case $p\in(0,1)$. For $p\in(0,1)$ it is known that the following variant of the reverse triangle inequality holds for all $x, y \in \R^d$ \citep[Exercise~24, Chapter~3]{Rudin1987}
	\[	\left\vert \left\Vert x \right\Vert_p^p - \left\Vert y \right\Vert_p^p \right\vert \leq \left\Vert x - y \right\Vert_p^p.	\]
Using this result, for $p\in(0,1)$ we can proceed analogously to \eqref{term 1} and write
	\begin{equation}	\label{term 1b}
	\begin{aligned}
	\left\vert \frac{\left\Vert \widetilde{u}(x) \right\Vert_p^p}{\left\Vert u \right\Vert_p^p} - 1 \right\vert \leq \frac{\left\Vert u - \widetilde{u}(x) \right\Vert_p^p}{\left\Vert u \right\Vert_p^p} \leq d^{1-p/2} \left\Vert u - \widetilde{u}(x) \right\Vert_2^p,
	\end{aligned}
	\end{equation}
as by \eqref{l_p relations} we have that
	\begin{equation}	\label{u1 inequality}
	\left\Vert u \right\Vert_p^p \geq \left\Vert u \right\Vert_2^p = 1	
	\end{equation}
that holds true for any $p \in (0,2]$. By \eqref{l_p relations} again and the previous inequality we know that
	\begin{equation}	\label{u inequality}
	\frac{\left\Vert \widetilde{u}(x) \right\Vert_p^p}{\left\Vert u \right\Vert_p^p} \leq d^{1-p/2}.	
	\end{equation}
For $0 \leq x \leq d^{1-p/2}$ and $p \in (0,1)$ we have by the mean value theorem that
	\[	\left\vert x^{1/p} - 1 \right\vert \leq \frac{d^{(1-p/2)(1/p-1)}}{p} \left\vert x - 1 \right\vert.	\]
It means that by \eqref{term 1b} and \eqref{u inequality} we can write
	\begin{equation}	\label{term 1c}
	\left\vert \frac{\left\Vert \widetilde{u}(x) \right\Vert_p}{\left\Vert u \right\Vert_p} - 1 \right\vert \leq \frac{d^{(1-p/2)(1/p-1)}}{p} \left\vert \frac{\left\Vert \widetilde{u}(x) \right\Vert_p^p}{\left\Vert u \right\Vert_p^p} - 1 \right\vert \leq \frac{d^{1/p-1/2}}{p} \left\Vert u - \widetilde{u}(x) \right\Vert_2^p,	
	\end{equation}
which gives the final inequality for the first term of \eqref{r ratio} also for $p \in (0,1)$.

The second term in \eqref{r ratio} can be treated at once for all $p \in (0,2]$. One can write
	\[	\frac{\left\langle x, u \right\rangle}{\left\langle x, \widetilde{u}(x) \right\rangle} = 1 + \frac{\left\langle x, u - \widetilde{u}(x) \right\rangle}{\left\langle x, \widetilde{u}(x) \right\rangle}	\]
and bound the second summand using \eqref{q norm argument}, \eqref{l_p relations}, and \eqref{u1 inequality} by
	\begin{equation}	\label{term 2}
	\left\vert \frac{\left\langle x, u - \widetilde{u}(x) \right\rangle}{\left\langle x, \widetilde{u}(x) \right\rangle} \right\vert \leq \frac{\left\Vert x \right\Vert_2 \left\Vert u - \widetilde{u}(x) \right\Vert_2}{\left\vert \left\langle x, \widetilde{u}(x) \right\rangle \right\vert} \leq \frac{\left\Vert x \right\Vert_2 \left\Vert u - \widetilde{u}(x) \right\Vert_2}{\left\Vert x \right\Vert_q \left\Vert \widetilde{u}(x) \right\Vert_p} \leq d^{1/2-1/q}\left\Vert u - \widetilde{u}(x) \right\Vert_2.	
	\end{equation}
	
Joining \eqref{term 1} and \eqref{term 2} together, for $p \in [1,2]$ we obtain for any $x \in \R^d$
	\begin{equation}	\label{final bound 1}
	\begin{aligned}
	\left\vert r(x,u) - 1 \right\vert & = \left\vert \frac{\left\Vert \widetilde{u}(x) \right\Vert_p}{\left\Vert u \right\Vert_p} \left( \frac{\left\langle x, u \right\rangle}{\left\langle x, \widetilde{u}(x) \right\rangle} - 1 + 1\right) - 1 \right\vert \\
	&	\leq \frac{\left\Vert \widetilde{u}(x) \right\Vert_p}{\left\Vert u \right\Vert_p} \left\vert \frac{\left\langle x, u \right\rangle}{\left\langle x, \widetilde{u}(x) \right\rangle} - 1 \right\vert + \left\vert \frac{\left\Vert \widetilde{u}(x) \right\Vert_p}{\left\Vert u \right\Vert_p} - 1 \right\vert \\
		& \leq d^{1/p-1/2} \left( d^{1/2-1/q} + 1\right) \left\Vert u - \widetilde{u}(x) \right\Vert_2\\
		& \leq 2d^{1/p-1/2} \left( d^{1/2-1/q} + 1\right)\sin\left(\left\Vert u - \widetilde{u}(x) \right\Vert_g /2\right)
	\end{aligned}
	\end{equation}
and the desired \eqref{s uniform convergence} is verified. For $p \in (0,1)$ we combine \eqref{term 1c} and \eqref{term 2} and get
	\begin{equation}	\label{final bound 2}
	\begin{aligned}
	\left\vert r(x,u) - 1 \right\vert & \leq \frac{\left\Vert \widetilde{u}(x) \right\Vert_p}{\left\Vert u \right\Vert_p} \left\vert \frac{\left\langle x, u \right\rangle}{\left\langle x, \widetilde{u}(x) \right\rangle} - 1 \right\vert + \left\vert \frac{\left\Vert \widetilde{u}(x) \right\Vert_p}{\left\Vert u \right\Vert_p} - 1 \right\vert \\
		& \leq d^{1/p-1/q}\left\Vert u - \widetilde{u}(x) \right\Vert_2 + \frac{d^{1/p-1/2}}{p} \left\Vert u - \widetilde{u}(x) \right\Vert_2^p\\
		& \leq 2d^{1/p-1/q}\sin\left(\left\Vert u - \widetilde{u}(x) \right\Vert_g /2\right) + 2^p\frac{d^{1/p-1/2}}{p} \sin^p\left(\left\Vert u - \widetilde{u}(x) \right\Vert_g /2\right).
	\end{aligned}
	\end{equation}

The statement about the rates of convergence follows by application of the proof technique from Theorem~\ref{theorem:modulus}, bounds \eqref{final bound 1} and \eqref{final bound 2}, and the relation \eqref{g 2 norm relation}.

\section{Proofs of the theoretical results: Projection depth}	\label{appendix:projection depth}

\subsection{Proof of Theorem~\ref{theorem:projection depth rate}}	\label{appendix:projection depth rate}

Take $M > \lambda_m$ and denote $B_M = \left\{ x \in \R^d \colon \left\Vert x \right\Vert < M \right\}$. Consider $x \in B_M$ first. For such $x$ we can write
	\begin{equation*}	
	\begin{aligned}
	\left\vert \varphi_x^P(u) - \varphi_x^P(v) \right\vert & \leq \left\vert \frac{\left\langle u, x \right\rangle - m_u}{s_u} - \frac{\left\langle u, x \right\rangle - m_u}{s_v} \right\vert + \left\vert \frac{\left\langle u, x \right\rangle - m_u}{s_v} - \frac{\left\langle v, x \right\rangle - m_v}{s_v} \right\vert \\
		& \leq \left( \left\vert \left\langle u, x \right\rangle \right\vert + \left\vert m_u \right\vert \right) \left\vert \frac{1}{s_u} - \frac{1}{s_v} \right\vert + \frac{1}{s_v} \left( \left\vert \left\langle u - v, x \right\rangle \right\vert + \left\vert m_v - m_u \right\vert \right) \\
		& \leq \left( \left\Vert x \right\Vert + \left\vert m_u \right\vert \right) \frac{1}{\lambda_{s,1}^2} \left\vert s_u - s_v \right\vert + \frac{1}{\lambda_{s,1}} \left( \left\Vert u - v \right\Vert \left\Vert x \right\Vert + \left\vert m_v - m_u \right\vert \right) \\
		& \leq \left( \left\Vert x \right\Vert + \lambda_m \right) \frac{1}{\lambda_{s,1}^2} C_s \left\Vert u - v \right\Vert + \frac{1}{\lambda_{s,1}} \left( \left\Vert u - v \right\Vert \left\Vert x \right\Vert + C_m \left\Vert u - v \right\Vert \right) \\
		& = \left\Vert u - v \right\Vert \left( \left\Vert x \right\Vert \left( \frac{C_s}{\lambda_{s,1}^2} + \frac{1}{\lambda_{s,1}} \right) + \frac{\lambda_m C_s}{\lambda_{s,1}^2} + \frac{C_m}{\lambda_{s,1}} \right) = \left\Vert u - v \right\Vert \left( \left\Vert x \right\Vert A + B \right) \\
		& \leq \left\Vert u - v \right\Vert \left( M A + B \right),
	\end{aligned}
	\end{equation*}
where the constants $A, B > 0$ are given by the last equality. From this it follows that
	\begin{equation}	\label{inner rate}
	\begin{aligned}
	\Delta_n^P\left(B_M\right) & \leq \delta_c\left( \left\Vert u - v \right\Vert \left( M A + B \right) \right) = \delta_c\left( 2\sin\left(\left\Vert u - v \right\Vert_g/2\right) \left( M A + B \right) \right) \\
	& \leq \delta_c\left( \left\Vert u - v \right\Vert_g \left( M A + B \right) \right),
	\end{aligned}	
	\end{equation}
where $\delta_c$ is the minimal modulus of continuity of $c$.

Now suppose that $x \notin B_M$. In what follows we bound, with $\widetilde{u}^P(x)$ abbreviated to $\widetilde{u}$ for simplicity, the expression
	\begin{equation}	\label{outer phi ratio}
	\left\vert \frac{\varphi_x^P(u)}{\varphi_x^P(\widetilde{u})} - 1 \right\vert \leq \left\vert \frac{s_{\widetilde{u}}}{s_u} \right\vert \left\vert \frac{ \left\langle x, u \right\rangle - m_u }{ \left\langle x, \widetilde{u} \right\rangle - m_{\widetilde{u}}} - 1 \right\vert + \left\vert \frac{s_{\widetilde{u}}}{s_u} - 1 \right\vert	
	\end{equation}
in terms of $\left\Vert \widetilde{u} - u \right\Vert_g$. First of all, note that $\left\vert s_{\widetilde{u}}/s_u \right\vert \leq \lambda_{s,2}/\lambda_{s,1}$. Secondly, by the definition of the maximizer $\widetilde{u} = \widetilde{u}^P(x)$, we have that
	\[	\varphi_x^P(\widetilde{u}) = \frac{ \left\langle x, \widetilde{u} \right\rangle - m_{\widetilde{u}} }{s_{\widetilde{u}}} \geq \frac{\left\langle x, x/\left\Vert x \right\Vert \right\rangle - m_{x/\left\Vert x \right\Vert}}{s_{x/\left\Vert x \right\Vert}} \geq \frac{\left\Vert x \right\Vert - \lambda_m}{\lambda_{s,2}} \geq \frac{M - \lambda_m}{\lambda_{s,2}} \geq 0,	\]
and thus
	\[	\left\vert \left\langle x, \widetilde{u} \right\rangle - m_{\widetilde{u}} \right\vert = \left\langle x, \widetilde{u} \right\rangle - m_{\widetilde{u}} \geq \frac{\lambda_{s,1}}{\lambda_{s,2}} \left(\left\Vert x \right\Vert - \lambda_m \right) \geq \frac{\lambda_{s,1}}{\lambda_{s,2}} \left(M - \lambda_m \right).	\] 
We can therefore write 
	\[	
	\begin{aligned}
	\left\vert \frac{ \left\langle x, u \right\rangle - m_u }{ \left\langle x, \widetilde{u} \right\rangle - m_{\widetilde{u}}} - 1 \right\vert & = \left\vert \frac{\left\langle x, u - \widetilde{u} \right\rangle + m_{\widetilde{u}} - m_u }{ \left\langle x, \widetilde{u} \right\rangle - m_{\widetilde{u}}} \right\vert \leq \frac{C_m \left\Vert \widetilde{u} - u \right\Vert_g + \left\Vert x \right\Vert \left\Vert \widetilde{u} - u \right\Vert}{\left\vert \left\langle x, \widetilde{u} \right\rangle - m_{\widetilde{u}} \right\vert} \\
	& \leq \frac{\lambda_{s,2}}{\lambda_{s,1}} \left( \frac{C_m \left\Vert \widetilde{u} - u \right\Vert_g}{M - \lambda_m} + \frac{\left\Vert x \right\Vert \left\Vert \widetilde{u} - u \right\Vert}{\left\Vert x \right\Vert - \lambda_m} \right) \\
	& = \frac{\lambda_{s,2}}{\lambda_{s,1}} \left( \frac{C_m \left\Vert \widetilde{u} - u \right\Vert_g}{M - \lambda_m} + \left\Vert \widetilde{u} - u \right\Vert + \frac{ \lambda_m \left\Vert \widetilde{u} - u \right\Vert}{\left\Vert x \right\Vert - \lambda_m} \right) \\
	& \le \frac{\lambda_{s,2}}{\lambda_{s,1}} \left( \frac{C_m \left\Vert \widetilde{u} - u \right\Vert_g}{M - \lambda_m} + \left\Vert \widetilde{u} - u \right\Vert + \frac{ \lambda_m \left\Vert \widetilde{u} - u \right\Vert}{M - \lambda_m} \right).
	\end{aligned}	\]

The second summand in \eqref{outer phi ratio} can be bounded by
	\[	\left\vert \frac{s_{\widetilde{u}}}{s_u} - 1 \right\vert = \left\vert \frac{s_{\widetilde{u}} - s_u}{s_u} \right\vert \leq \frac{C_s \left\Vert \widetilde{u} - u \right\Vert_g}{\lambda_{s,1}}.	\]
	
Altogether, we can write
	\[	\left\vert \frac{\varphi_x^P(u)}{\varphi_x^P(\widetilde{u})} - 1 \right\vert \leq \left(\frac{\lambda_{s,2}}{\lambda_{s,1}}\right)^2 \frac{C_m \left\Vert \widetilde{u} - u \right\Vert_g + M \left\Vert \widetilde{u} - u \right\Vert}{M - \lambda_m} + \frac{C_s \left\Vert \widetilde{u} - u \right\Vert_g}{\lambda_{s,1}}.	\]
This means that for $u \in \Sph$ such that $\left\Vert \widetilde{u} - u \right\Vert_g \leq \varepsilon$, i.e. $\left\Vert \widetilde{u} - u \right\Vert \leq 2\sin(\varepsilon/2) \leq \varepsilon$, we have
	\[	
	\begin{aligned}
	\left\vert \varphi_x^P(u) - \varphi_x^P(\widetilde{u})	\right\vert & \leq \varphi_x^P(\widetilde{u}) \left( \left(\frac{\lambda_{s,2}}{\lambda_{s,1}}\right)^2 \frac{\varepsilon \, C_m + 2M\sin(\varepsilon/2)}{M - \lambda_m} + \frac{\varepsilon \, C_s}{\lambda_{s,1}} \right) \\
	& \leq \varepsilon \, \varphi_x^P(\widetilde{u}) \left( \left(\frac{\lambda_{s,2}}{\lambda_{s,1}}\right)^2 \frac{ C_m + M }{M - \lambda_m} + \frac{ C_s}{\lambda_{s,1}} \right)
	\end{aligned}
	\]
or, because $\varphi_x^P(\widetilde{u}) \geq \varphi_x^P(u)$ for all $u \in \Sph$ by the definition of $\widetilde{u}$, 
	\[  
	\begin{aligned}
	\varphi_x^P(u) & \geq \varphi_x^P(\widetilde{u}) \left( 1 - \left(\frac{\lambda_{s,2}}{\lambda_{s,1}}\right)^2 \frac{\varepsilon \, C_m + 2M \sin(\varepsilon/2)}{M - \lambda_m} - \frac{\varepsilon \, C_s}{\lambda_{s,1}} \right) \\
	& \geq \varphi_x^P(\widetilde{u}) \left( 1 - \varepsilon \left( \left(\frac{\lambda_{s,2}}{\lambda_{s,1}}\right)^2 \frac{ C_m + M }{M - \lambda_m} + \frac{ C_s}{\lambda_{s,1}} \right) \right).
	\end{aligned}
	\]
For any $n \in 1, 2, \dots$ such that $\min_{i=1,\dots,n} \left\Vert U_i - \widetilde{u} \right\Vert_g \leq \varepsilon$ for all $x\notin B_M$ the last bounds result in 
	\[
	\begin{aligned}
	\PD(x;P) & = c\left( \varphi_x^P(\widetilde{u}) \right), \\
	\PD_n(x;P) & \leq c(\varphi_x^P(u)) \\
	& \leq c\left( \varphi_x^P(\widetilde{u}) \left( 1 - \left(\frac{\lambda_{s,2}}{\lambda_{s,1}}\right)^2 \frac{\varepsilon \, C_m + 2M \sin(\varepsilon/2)}{M - \lambda_m} - \frac{\varepsilon \, C_s}{\lambda_{s,1}} \right) \right) \\
	& \leq c\left( \varphi_x^P(\widetilde{u}) \left( 1 - \varepsilon \left( \left(\frac{\lambda_{s,2}}{\lambda_{s,1}}\right)^2 \frac{ C_m + M }{M - \lambda_m} + \frac{ C_s}{\lambda_{s,1}} \right) \right) \right).
	\end{aligned}
	\]
This holds true for any $x \notin B_M$. Consequently, if $\min_{i=1,\dots,n} \left\Vert U_i - \widetilde{u} \right\Vert_g \leq \varepsilon$, we can write
	\begin{equation}	\label{outer rate}
	\Delta_n^P\left(\R^d \setminus B_M\right) \leq c\left( \varphi_x^P(\widetilde{u}) \tau \right) - c(\varphi_x^P(\widetilde{u})) \leq \zeta(\tau_1) \leq \zeta(\tau_2)	
	\end{equation}
for 
	\[
	\begin{aligned}
	\tau_1 & = 1 - \left(\frac{\lambda_{s,2}}{\lambda_{s,1}}\right)^2 \frac{\varepsilon \, C_m + 2M \sin(\varepsilon/2)}{M - \lambda_m} - \frac{\varepsilon \, C_s}{\lambda_{s,1}}, \\
	\tau_2 & = 1 - \varepsilon \left( \left(\frac{\lambda_{s,2}}{\lambda_{s,1}}\right)^2 \frac{ C_m + M }{M - \lambda_m} + \frac{ C_s}{\lambda_{s,1}} \right).
	\end{aligned}
	\]

Finally, the rate of convergence of the approximations of the projection depth can be obtained as a combination of \eqref{inner rate} and \eqref{outer rate}
	\begin{equation}	\label{projection depth final rate}
	\begin{aligned}
	\Delta_n^P(\R^d) & \leq \inf_{M>\lambda_m} \max\left\{ \delta_c\left( 2\sin(\varepsilon/2) \left( M A + B \right) \right), \zeta(\tau_1)	\right\} \\
	& \leq \inf_{M>\lambda_m} \max\left\{ \delta_c\left( \varepsilon \left( M A + B \right) \right), \zeta(\tau_2)	\right\},
	\end{aligned}
	\end{equation}
where $\varepsilon$ is substituted by the 
maximal spacing $S_n$ introduced in the proof of Theorem~\ref{theorem:Lipschitz}. The final rates of convergence then follow by application of the same strategy as in the proof of Theorem~\ref{theorem:modulus}, where the result of \citet{Janson1987} is employed to devise an almost sure upper bound on the rate of convergence of the depth approximation. Note that to obtain explicit, almost sure bounds such as those from Theorems~\ref{theorem:Lipschitz} and~\ref{theorem:modulus} it is necessary to invert (an upper bound to) one of the functions on the right hand side of \eqref{projection depth final rate} considered as a function of $\varepsilon$.

\subsection{Proof of Theorem~\ref{theorem:projection depth spherical}}	\label{appendix:projection depth spherical}

Due to the translation invariance of $\PD$, the distribution $P$ can be assumed to be spherically symmetric about the origin. Then, it can be seen that $m_u = 0$ and $s_u = S$ for some constant $S \geq 0$. The latter fact follows because for such spherically symmetric distributions, all projections of $P$ onto lines passing through the origin have the same univariate distribution \citep{Serfling_symmetry}. If $S = 0$, the assertion of the theorem is trivially satisfied, as in this case $\PD(x;P) = \PD_n(x;P) = 0$ for all $x \in \R^d$. For $S>0$ for any $x \in \R^d$
	\begin{equation*}	
	\PD(x;P) = c\left( \sup_{u \in \Sph} \frac{\left\langle x, u \right\rangle}{S} \right) = c\left( \left\Vert x \right\Vert/S \right),	
	\end{equation*}
and the depth is realised in $\widetilde{u}^P(x) = x/\left\Vert x \right\Vert$. Likewise, for function $\varphi_x^P$ from the approximating projection depth it holds true that
	\[	\varphi_x^P(u) = \frac{\left\langle x, u \right\rangle}{S}.	\]
Similarly as in the proof of Theorem~\ref{theorem:projection depth rate} let us bound the expression bounding \eqref{outer phi ratio}. Here we can for $x \ne 0$ write
	\[	\left\vert \frac{\varphi_x^P(u)}{\varphi_x^P(\widetilde{u}^P(x))} - 1 \right\vert = \left\vert \frac{\left\langle x, u \right\rangle}{\left\langle x, \widetilde{u}^P(x) \right\rangle} - 1 \right\vert = \frac{\left\vert \left\langle x, u - \widetilde{u}^P(x) \right\rangle \right\vert}{\left\Vert x \right\Vert} \leq \left\Vert u - \widetilde{u}^P(x) \right\Vert \leq \left\Vert u - \widetilde{u}^P(x) \right\Vert_g.	\]
Thus, in the same way as in \eqref{outer rate} in the proof of Theorem~\ref{theorem:projection depth rate} we obtain that under the maximal spacing condition $\min_{i=1,\dots,n} \left\Vert U_i - \widetilde{u} \right\Vert_g \leq \varepsilon$ we can write
	\[	\Delta_n^P\left(\R^d\right) \leq c\left( \varphi_x^P(\widetilde{u}^P(x)) \tau \right) - c(\varphi_x^P(\widetilde{u}^P(x))) \leq \zeta(\tau_1) \leq \zeta(\tau_2),	\]
where
	\[
	\begin{aligned}
	\tau_1 & = 1 - 2\sin(\varepsilon/2), \\
	\tau_2 & = 1 - \varepsilon.
	\end{aligned}
	\]
For the final rates we apply the same technique as in the proof of Theorem~\ref{theorem:modulus}; note that in the statement of the theorem only the second, simpler rate is stated.

The formula for $\zeta$ in \eqref{zeta function} can be obtained by direct maximization in the expression for the multiplicative modulus of continuity.

\subsection*{Acknowledgement}
Stanislav Nagy was supported by the grant 19-16097Y of the Czech Science Foundation, and by the PRIMUS/17/SCI/3 project of Charles University.


\def\cprime{$'$} \def\polhk#1{\setbox0=\hbox{#1}{\ooalign{\hidewidth
  \lower1.5ex\hbox{`}\hidewidth\crcr\unhbox0}}}


\begin{thebibliography}{}

\bibitem[Aaron et~al., 2017]{Aaron_etal2017}
Aaron, C., Cholaquidis, A., and Fraiman, R. (2017).
\newblock A generalization of the maximal-spacings in several dimensions and a
  convexity test.
\newblock {\em Extremes}, 20(3):605--634.

\bibitem[Bogi\'{c}evi\'{c} and Merkle, 2018]{Bogicevic_Merkle2018}
Bogi\'{c}evi\'{c}, M. and Merkle, M. (2018).
\newblock Approximate calculation of {T}ukey's depth and median with
  high-dimensional data.
\newblock {\em Yugosl. J. Oper. Res.}, 28(4):475--499.

\bibitem[Burr and Fabrizio, 2017]{Burr_Fabrizio2017}
Burr, M.~A. and Fabrizio, R.~J. (2017).
\newblock Uniform convergence rates for halfspace depth.
\newblock {\em Statist. Probab. Lett.}, 124:33--40.

\bibitem[Chen et~al., 2013]{Chen_etal2013}
Chen, D., Morin, P., and Wagner, U. (2013).
\newblock Absolute approximation of {T}ukey depth: theory and experiments.
\newblock {\em Comput. Geom.}, 46(5):566--573.

\bibitem[Chen and Tyler, 2004]{Chen_Tyler2004}
Chen, Z. and Tyler, D.~E. (2004).
\newblock On the behavior of {T}ukey's depth and median under symmetric stable
  distributions.
\newblock {\em J. Statist. Plann. Inference}, 122(1-2):111--124.

\bibitem[Cuesta-Albertos and Nieto-Reyes, 2008]{Cuesta_Nieto2008}
Cuesta-Albertos, J.~A. and Nieto-Reyes, A. (2008).
\newblock The random {T}ukey depth.
\newblock {\em Comput. Statist. Data Anal.}, 52(11):4979--4988.

\bibitem[DeVore and Lorentz, 1993]{deVore_Lorentz1993}
DeVore, R.~A. and Lorentz, G.~G. (1993).
\newblock {\em Constructive approximation}, volume 303 of {\em Grundlehren der
  Mathematischen Wissenschaften [Fundamental Principles of Mathematical
  Sciences]}.
\newblock Springer-Verlag, Berlin.

\bibitem[Devroye, 1981]{Devroye1981}
Devroye, L. (1981).
\newblock Laws of the iterated logarithm for order statistics of uniform
  spacings.
\newblock {\em Ann. Probab.}, 9(5):860--867.

\bibitem[Diaconis and Freedman, 1984]{Diaconis_Freedman1984}
Diaconis, P. and Freedman, D. (1984).
\newblock Asymptotics of graphical projection pursuit.
\newblock {\em Ann. Statist.}, 12(3):793--815.

\bibitem[Donoho, 1982]{Donoho1982}
Donoho, D.~L. (1982).
\newblock Breakdown properties of multivariate location estimators.
\newblock Qualifying paper, Harvard University.

\bibitem[Donoho and Gasko, 1992]{Donoho_Gasko1992}
Donoho, D.~L. and Gasko, M. (1992).
\newblock Breakdown properties of location estimates based on halfspace depth
  and projected outlyingness.
\newblock {\em Ann. Statist.}, 20(4):1803--1827.

\bibitem[Dudley, 2002]{Dudley2002}
Dudley, R.~M. (2002).
\newblock {\em Real analysis and probability}, volume~74 of {\em Cambridge
  Studies in Advanced Mathematics}.
\newblock Cambridge University Press, Cambridge.
\newblock Revised reprint of the 1989 original.

\bibitem[Dyckerhoff, 2004]{Dyckerhoff2004}
Dyckerhoff, R. (2004).
\newblock Data depths satisfying the projection property.
\newblock {\em Allg. Stat. Arch.}, 88(2):163--190.

\bibitem[Embrechts and Hofert, 2013]{Embrechts_Hofert2013}
Embrechts, P. and Hofert, M. (2013).
\newblock A note on generalized inverses.
\newblock {\em Math. Methods Oper. Res.}, 77(3):423--432.

\bibitem[Fang et~al., 1990]{Fang_etal1990}
Fang, K.~T., Kotz, S., and Ng, K.~W. (1990).
\newblock {\em Symmetric multivariate and related distributions}, volume~36 of
  {\em Monographs on Statistics and Applied Probability}.
\newblock Chapman and Hall, Ltd., London.

\bibitem[Horn and Johnson, 1994]{Horn_Johnson1994}
Horn, R.~A. and Johnson, C.~R. (1994).
\newblock {\em Topics in matrix analysis}.
\newblock Cambridge University Press, Cambridge.
\newblock Corrected reprint of the 1991 original.

\bibitem[Janson, 1986]{Janson1986}
Janson, S. (1986).
\newblock Random coverings in several dimensions.
\newblock {\em Acta Math.}, 156(1-2):83--118.

\bibitem[Janson, 1987]{Janson1987}
Janson, S. (1987).
\newblock Maximal spacings in several dimensions.
\newblock {\em Ann. Probab.}, 15(1):274--280.

\bibitem[Johnson and Preparata, 1978]{Johnson_Prepata1978}
Johnson, D.~S. and Preparata, F.~P. (1978).
\newblock The densest hemisphere problem.
\newblock {\em Theoret. Comput. Sci.}, 6(1):93--107.

\bibitem[Kuelbs and Dudley, 1980]{Kuelbs_Dudley1980}
Kuelbs, J. and Dudley, R.~M. (1980).
\newblock Log log laws for empirical measures.
\newblock {\em Ann. Probab.}, 8(3):405--418.

\bibitem[Lin, 1972]{Lin1972}
Lin, P.~E. (1972).
\newblock Some characterizations of the multivariate {$t$} distribution.
\newblock {\em J. Multivariate Anal.}, 2:339--344.

\bibitem[Liu et~al., 1999]{Liu_etal1999}
Liu, R.~Y., Parelius, J.~M., and Singh, K. (1999).
\newblock Multivariate analysis by data depth: descriptive statistics, graphics
  and inference.
\newblock {\em Ann. Statist.}, 27(3):783--858.

\bibitem[Liu and Singh, 1993]{Liu_Singh1993}
Liu, R.~Y. and Singh, K. (1993).
\newblock A quality index based on data depth and multivariate rank tests.
\newblock {\em J. Amer. Statist. Assoc.}, 88(421):252--260.

\bibitem[Liu and Zuo, 2014]{Liu_Zuo2014}
Liu, X. and Zuo, Y. (2014).
\newblock Computing projection depth and its associated estimators.
\newblock {\em Stat. Comput.}, 24(1):51--63.

\bibitem[Mass{\'e}, 2004]{Masse2004}
Mass{\'e}, J.-C. (2004).
\newblock Asymptotics for the {T}ukey depth process, with an application to a
  multivariate trimmed mean.
\newblock {\em Bernoulli}, 10(3):397--419.

\bibitem[Mass{\'e} and Theodorescu, 1994]{Masse_Theodorescu1994}
Mass{\'e}, J.-C. and Theodorescu, R. (1994).
\newblock Halfplane trimming for bivariate distributions.
\newblock {\em J. Multivariate Anal.}, 48(2):188--202.

\bibitem[Mosler, 2002]{Mosler2002}
Mosler, K. (2002).
\newblock {\em Multivariate dispersion, central regions and depth: The lift
  zonoid approach}, volume 165 of {\em Lecture Notes in Statistics}.
\newblock Springer-Verlag, Berlin.

\bibitem[Mozharovskyi et~al., 2015]{Mozharovskyi_etal2015}
Mozharovskyi, P., Mosler, K., and Lange, T. (2015).
\newblock Classifying real-world data with the {$DD\alpha$}-procedure.
\newblock {\em Adv. Data Anal. Classif.}, 9(3):287--314.

\bibitem[Rousseeuw and Ruts, 1999]{Rousseeuw_Ruts1999}
Rousseeuw, P.~J. and Ruts, I. (1999).
\newblock The depth function of a population distribution.
\newblock {\em Metrika}, 49(3):213--244.

\bibitem[Rudin, 1987]{Rudin1987}
Rudin, W. (1987).
\newblock {\em Real and complex analysis}.
\newblock McGraw-Hill Book Co., New York, third edition.

\bibitem[Serfling, 2006]{Serfling_symmetry}
Serfling, R. (2006).
\newblock Multivariate symmetry and asymmetry.
\newblock {\em Encyclopedia of Statistical Sciences, Second Edition},
  8:5338--5345.

\bibitem[Shao and Zuo, 2012]{Shao_Zuo2012}
Shao, W. and Zuo, Y. (2012).
\newblock Simulated annealing for higher dimensional projection depth.
\newblock {\em Comput. Statist. Data Anal.}, 56(12):4026--4036.

\bibitem[Shao and Zuo, 2019]{Shao_Zuo2019}
Shao, W. and Zuo, Y. (2019).
\newblock Computing the halfspace depth with multiple try algorithm and
  simulated annealing algorithm.
\newblock {\em Computational Statistics}.
\newblock To appear.

\bibitem[Stahel, 1981]{Stahel1981}
Stahel, W.~A. (1981).
\newblock Robuste {S}ch\"atzungen: {I}nfinitesimale {O}ptimalit\"at und
  {S}ch\"atzungen von {K}ovarianzmatrizen.
\newblock Ph.D. thesis, ETH Z\"urich.

\bibitem[Tukey, 1975]{Tukey1975}
Tukey, J.~W. (1975).
\newblock Mathematics and the picturing of data.
\newblock In {\em Proceedings of the {I}nternational {C}ongress of
  {M}athematicians ({V}ancouver, {B}. {C}., 1974), {V}ol. 2}, pages 523--531.
  Canad. Math. Congress, Montreal, Que.

\bibitem[Yukich, 1990]{Yukich1990}
Yukich, J.~E. (1990).
\newblock The law of the iterated logarithm for empirical processes.
\newblock In {\em Probability in {B}anach spaces 6 ({S}andbjerg, 1986)},
  volume~20 of {\em Progr. Probab.}, pages 265--282. Birkh\"{a}user Boston,
  Boston, MA.

\bibitem[Zolotarev, 1986]{Zolotarev1986}
Zolotarev, V.~M. (1986).
\newblock {\em One-dimensional stable distributions}, volume~65 of {\em
  Translations of Mathematical Monographs}.
\newblock American Mathematical Society, Providence, RI.
\newblock Translated from the Russian by H. H. McFaden, Translation edited by
  Ben Silver.

\bibitem[Zuo, 2003]{Zuo2003}
Zuo, Y. (2003).
\newblock Projection-based depth functions and associated medians.
\newblock {\em Ann. Statist.}, 31(5):1460--1490.

\bibitem[Zuo and Lai, 2011]{Zuo_Lai2011}
Zuo, Y. and Lai, S. (2011).
\newblock Exact computation of bivariate projection depth and the
  {S}tahel-{D}onoho estimator.
\newblock {\em Comput. Statist. Data Anal.}, 55(3):1173--1179.

\bibitem[Zuo and Serfling, 2000]{Zuo_Serfling2000}
Zuo, Y. and Serfling, R. (2000).
\newblock General notions of statistical depth function.
\newblock {\em Ann. Statist.}, 28(2):461--482.

\end{thebibliography}
\end{document}